\documentclass[12pt]{amsart}
\usepackage{a4wide}
\usepackage[centertags,leqno]{amsmath}
\usepackage{amssymb,amsfonts}
\usepackage{color}
\usepackage{graphicx}
\usepackage[colorlinks=true]{hyperref}

\newtheorem{theorem}{Theorem}[section]

\newtheorem{corollary}[theorem]{Corollary}

\newtheorem{lemma}[theorem]{Lemma}
\newtheorem{proposition}[theorem]{Proposition}
\newtheorem{remark}[theorem]{Remark}

\definecolor{violet}{rgb}{0.5,0,0.5}
\definecolor{orange}{cmyk}{0,0.3,0.7,0}

\newtheorem{definition}[theorem]{Definition}

\numberwithin{equation}{section}

\renewcommand{\qed}{\rule{2mm}{2mm}}

\newcommand{\eqdef}{\stackrel{{\mathrm {def}}}{=}}
\newcommand{\eps}{\varepsilon}
\renewcommand{\colon}{:\,}
\renewcommand{\proof}{{\em Proof. }}

\newcommand{\RR}{\mathbb{R}}

\newcommand{\Div}{\mathop{\rm {div}}}

\newtheorem{Theorem}{Theorem}[section]
\newtheorem{Proposition}{Proposition}[section]
\newtheorem{Lemma}{Lemma}[section]
\newtheorem{Corollary}{Corollary}[section]
\newtheorem{Definition}{Definition}[section]

\newcommand{\bTheorem}[1]{
\begin{Theorem} \label{T#1} }
\newcommand{\eT}{\end{Theorem}}

\newcommand{\bProposition}[1]{
\begin{Proposition} \label{P#1}}
\newcommand{\eP}{\end{Proposition}}

\newcommand{\bLemma}[1]{
\begin{Lemma} \label{L#1} }
\newcommand{\eL}{\end{Lemma}}

\newcommand{\bCorollary}[1]{
\begin{Corollary} \label{C#1} }
\newcommand{\eC}{\end{Corollary}}

\newcommand{\bDefinition}[1]{
\begin{Definition} \label{D#1} }
\newcommand{\eD}{\end{Definition}}

\newcommand{\bFormula}[1]{
\begin{equation} \label{#1}}
\newcommand{\eF}{\end{equation}}

\newcommand{\bProof}{{\bf Proof: }}

\newcommand{\dx}{{\rm d} {x}}

\definecolor{Cgrey}{rgb}{0.85,0.85,0.85}

\newcommand\Cbox[2]{%
    \newbox\contentbox%
    \newbox\bkgdbox%
    \setbox\contentbox\hbox to \hsize{%
        \vtop{
            \kern\columnsep
            \hbox to \hsize{%
                \kern\columnsep%
                \advance\hsize by -2\columnsep%
                \setlength{\textwidth}{\hsize}%
                \vbox{
                    \parskip=\baselineskip
                    \parindent=0bp
                    #2
                }%
                \kern\columnsep%
            }%
            \kern\columnsep%
        }%
    }%
    \setbox\bkgdbox\vbox{
        \color{#1}
        \hrule width  \wd\contentbox %
               height \ht\contentbox %
               depth  \dp\contentbox
        \color{black}
    }%
    \wd\bkgdbox=0bp%
    \vbox{\hbox to \hsize{\box\bkgdbox\box\contentbox}}%
    \vskip\baselineskip%
}

\begin{document}

\title[ Front propagation in parabolic equations ]
{Front propagation in non\-linear parabolic equations}

\author{Eduard Feireisl}
\thanks{E.F. acknowledges support from
        the Czech Science Foundation (GA\v{C}R) under
        the project\hfill\break
        {P}201-13-00522{S} in the framework of RVO: 67985840.}
\address{
  Eduard Feireisl,
  Institute of Mathematics,
  Academy of Sciences of\hfil\break the Czech Republic,
  \v{Z}itn\'a~25,
  CZ - 115~67 Praha~1, Czech Republic
}
\email{feireisl@math.cas.cz}

\author{Danielle Hilhorst}
\address{
  Danielle Hilhorst,
  CNRS et Laboratoire de Math\'ematiques,
  Universit\'e de Paris\--Sud, B\^atiment 425,
  FR - 91405 Orsay Cedex, France
}
\email{Danielle.Hilhorst@math.u-psud.fr}

\author{Hana Petzeltov\'a}
\address{
  Hana Petzeltov\'a,
  Institute of Mathematics,
  Academy of Sciences of\hfil\break the Czech Republic,
  \v{Z}itn\'a~25,
  CZ - 115~67 Praha~1, Czech Republic
}
\email{petzelt@math.cas.cz}

\author{Peter Tak\'a\v{c}}
\thanks{P.T. was partially supported by the grant TA~213/15-1 from
        the German Research Foundation (D.F.G.) of Germany.}
\address{
  Peter Tak\'a\v{c},
  Institut f\"ur Mathematik,
  Universit\"at Rostock,
  Ulmenstra{\ss}e~69,\hfil\break Haus~3,
  D-18055 Rostock, Germany
}
\email{peter.takac@uni-rostock.de}

\maketitle

\vspace{0.5cm}


\baselineskip=12pt
\noindent
\begingroup\footnotesize
{\bf {\sc Abstract.}}
We study existence and stability of travelling waves for
non\-linear convection diffusion equations in the 1-D Euclidean space.
The diffusion coefficient depends on the gradient
in analogy with the $p$-Laplacian and may be degenerate.
Unconditional stability is established with respect to
initial data perturbations in $L^1(\mathbb{R})$.
\endgroup

\vfill
\par\vspace*{0.2cm}
\noindent
\begin{tabular}{ll}
{\bf Running head:}
& Poincar\'e inequality and P.-S.\ condition\\
\end{tabular}

\par\vspace*{0.2cm}
\noindent
\begin{tabular}{ll}
{\bf Keywords:}
& Poincar\'e inequality, Palais\--Smale condition,\\
& $p$-Laplacian, energy functional\\
\end{tabular}

\par\vspace*{0.2cm}
\noindent
\begin{tabular}{ll}
{\bf 2000 Mathematics Subject Classification:}
& Primary   35J20, 35B45;\\
& Secondary 35P30, 46E35 \\
\end{tabular}

\newpage

\baselineskip=14pt
\section{Introduction}
\label{s:Intro}

The purpose of this article is to investigate
the propagation of very simple {\it travelling waves\/}
in a reaction\--diffusion model.
The model is the favorite {\it Fisher\--KPP equation\/}
(or {\it Fisher\--Kolmogorov equation\/})
derived by {\sc R.~A.\ Fisher} \cite{Fisher} in $1937$
and first mathematically analyzed by
{\sc A.\ Kolmogorov}, {\sc I.\ Petrovski}, and {\sc N.\ Piscounov}
\cite{KPP} in the same year.
However, these original works (\cite{Fisher, KPP})
consider solely {\it linear\/} diffusion and
(sufficiently) {\it smooth\/} reaction.
In our present work, we allow for both,
a {\it nonlinear\/} diffusion operator and
a {\it nonsmooth\/} reaction function.
More precisely, we study the {\it\bfseries interaction\/}
between the (nonlinear) diffusion and the (nonsmooth) reaction;
in paticular, their influence on the formation and the shape of
a {\it travelling wave\/} connecting
two stable (spatially constant) steady states.

We consider the following non\-linear evolutionary problem
for an unknown function $u = u(x,t)$,
\begin{equation}
\label{i1}
  \partial_t u = \Div\left( \partial\Phi (\nabla_x u) \right) + f(u) \,,
  \quad x\in \mathbb{R}^N \,,\ t > 0 \,,
\end{equation}
supplemented by the initial condition
\begin{equation}
\label{i2}
  u(\,\cdot\,,0) = u_0 \quad\mbox{ in }\, \mathbb{R}^N \,.
\end{equation}
Here,
$\Phi\colon \mathbb{R}^N\to \mathbb{R}$,
$f\colon \mathbb{R}\to \mathbb{R}$, and
$u_0\colon \mathbb{R}^N\to \mathbb{R}$ are given data as specified below.
Roughly speaking, we assume that
$\Phi$ is a continuously differentiable, convex functional on
$\mathbb{R}^N$ with the Fr\'echet derivative
$\partial\Phi\colon \mathbb{R}^N\to \mathbb{R}^N$, such that
$\Phi$ is also radially symmetric of class
$C^2\left( \mathbb{R}^N\setminus \{\mathbf{0}\} \right)$,
its Hessian matrix
$\partial^2\Phi(Z)\in \mathbb{R}^{N\times N}$
is positively definite at every point
$Z\in \mathbb{R}^N\setminus \{\mathbf{0}\}$, and
\begin{equation*}
  | \partial^2\Phi(Z) |\cdot |Z|\to 0 \quad\mbox{ as }\, |Z|\to 0 \,.
\end{equation*}
The non\-linear reaction function
$f\colon \mathbb{R}\to \mathbb{R}$ is of the {\em\bfseries KPP\--type\/}
({\sc Kolmogorov}\--{\sc Pet\-rov\-ski}\--{\sc Piscounov} \cite{KPP});
specifically,
\begin{equation}
\label{i3}
  f\in C(\mathbb{R}) \,,\ f(-1) = f(\mu) = f(1) = 0 \,,\
  f < 0 \,\mbox{ in }\, (-1,\mu) \,,\
  f > 0 \,\mbox{ in }\, (\mu, 1) \,.
\end{equation}
Moreover, we assume that its integral
\begin{equation*}
  F(r)\eqdef \int_{-1}^r f(s) \,\mathrm{d}s \,,\quad -1\leq r\leq 1 \,,
\end{equation*}
satisfies
\begin{equation}
\label{int:f(u)}
  F(1) - F(r)
  = \int_{r}^{1} f(s) \,\mathrm{d}s > 0
    \quad\mbox{ whenever }\, -1 < r < 1 \,.
\end{equation}

Taking the initial data
$u_0\colon \mathbb{R}^N\to \mathbb{R}$
valued in the interval $[-1,1]$
between the extremal zeros ($= \mp 1$) of $f$, i.e., $-1\leq u_0\leq 1$,
we are interested in the long\--time behavior of solutions to
problem \eqref{i1}, \eqref{i2}; in particular,
in propagation of fronts separating the areas where
$u$ approaches the limit values $\pm 1$, respectively.

A currently standard approach consists of introducing
the (hyperbolic) change of variables
$t\approx \frac{t}{\eps}$ and $x\approx \frac{x}{\eps}$
which leads us to the scaled problem
\begin{equation}
\label{i4}
  \partial_t u_{\eps}
  = \Div \left( \partial \Phi (\eps\nabla_x u_{\eps} )  \right)
  + \frac{1}{\eps} f(u_{\eps})\,,
  \quad x\in \mathbb{R}^N \,,\ t > 0 \,,
\end{equation}
for the unknown function
\begin{math}
  u_{\eps}(x,t) \eqdef
  u\left( \frac{x}{\eps} ,\, \frac{t}{\eps} \right) ,
\end{math}
supplemented by the initial data
\begin{equation}
\label{i5}
  u_{\eps}(\,\cdot\,,0) = u_{\eps,0} \,\mbox{ in }\, \mathbb{R}^N \,,
\end{equation}
see the survey by {\sc P.~E.\ Souganidis} \cite{Sougan}.
A prototype example of \eqref{i4} is the equation involving
the $p$-{\em\bfseries Laplace operator\/} ($p > 1$),
\begin{equation}
\label{i6}
  \partial_t u_{\eps}
  = {\eps}^{p-1} \Div
    \left( |\nabla_x u_{\eps}|^{p-2} \nabla_x u_{\eps} \right)
  + \frac{1}{\eps} f(u_{\eps}) \,,
  \quad x\in \mathbb{R}^N \,,\ t > 0 \,.
\end{equation}

Our aim is to examine the behavior of solutions $u_{\eps}$ of
problem \eqref{i4}, \eqref{i5} as $\eps\to 0+$.
In particular, we extend the results of
{\sc Zhao} and {\sc Yi} \cite{ZhYi} for problem~\eqref{i6} with
$p > 2$ (the ``degenerate case'' of slow diffusion)
to the ``singular case'' $1 < p < 2$ of fast diffusion.

\subsection{Travelling waves}

The asymptotic behavior of solutions $u_{\eps}$ to problem \eqref{i6}
in the singular limit $\eps\to 0+$ is well\--understood in
the non\-degenerate case (slow diffusion case)
$p=2$ ($p > 2$) and also for
the porous media type elliptic operator $\Delta u^m$, $m > 1$, see
{\sc Aronson} and {\sc Weinberger} \cite{AronWein},
{\sc Chen} \cite{Chen-XF},
{\sc Feireisl} \cite{EF14},
{\sc Fife} and {\sc McLeod} \cite{FM}, and
{\sc Zhao} and {\sc Yi} \cite{ZhYi}.

To begin, we decompose $\mathbb{R}^N$ into
the closures of the following two regions:
\begin{equation}
\label{i7}
\begin{aligned}
  G_{-}\eqdef
    \left\{ x\in \mathbb{R}^N \colon
    \,\mbox{ there exists a neighborhood $U(x)\subset \mathbb{R}^N$
             of $x$ such that }\,
    \right.
\\
    \left.
  \limsup_{\eps\to 0+} \left( \sup_{y\in U(x)} u_{\eps,0}(y)\right) < \mu
    \right\} \,,
\end{aligned}
\end{equation}
and
\begin{equation}
\label{i8}
\begin{aligned}
  G_{+}\eqdef
    \left\{ x\in \mathbb{R}^N \colon
    \,\mbox{ there exists a neighborhood $U(x)\subset \mathbb{R}^N$
             of $x$ such that }\,
    \right.
\\
    \left.
  \limsup_{\eps\to 0+} \left( \inf_{y\in U(x)} u_{\eps,0}(y)\right) > \mu
    \right\} \,,
\end{aligned}
\end{equation}
where we assume that
\begin{equation}
\label{i9}
  \overline{G}_- \cup \overline{G}_+ = \mathbb{R}^N \,,
  \quad\mbox{ and set }\quad
  \Gamma = \overline{G}_- \cap\overline{G}_+ \,.
\end{equation}

Similarly to
{\sc Barles}, {\sc Bronsard}, and {\sc Souganidis} \cite{BBS},
we expect that
\begin{align}
\label{i10}
  u_{\eps}\to -1 &\,\mbox{ uniformly in compact subsets of }\,
  \left\{ (x,t)\in \mathbb{R^N}\times \RR_+\colon
          \mathrm{dist}[x, \Gamma] > c t
  \right\} \,,
\\
\label{i11}
  u_{\eps}\to 1 &\,\mbox{ uniformly in compact subsets of }\,
  \left\{ (x,t)\in \mathbb{R^N}\times \RR_+\colon
          \mathrm{dist}[x, \Gamma] < c t
  \right\} \,,
\end{align}
where $\RR_+\eqdef [0,\infty)$ and
``$\mathrm{dist}$'' stands for the {\em signed distance\/},
\begin{equation}
\label{DIST}
  \mathrm{dist}[x, \Gamma] =
\left\{
  \begin{aligned}
  \inf\{ |x-y|\colon y\in \Gamma\}
&   \quad\mbox{ for }\, x\in G_{-} \,,
\\
  {}- \inf\{ |x-y|\colon y\in \Gamma\}
&   \quad\mbox{ for }\, x\in G_{+} \,,
  \end{aligned}
\right.
\end{equation}
where $c$ is the \emph{front speed} that can be determined as
the speed of propagation of the traveling waves for
the associated $1$D problem.

Setting $u(x,t) = q(x-ct)$ in eq.~\eqref{i1} we have
$\partial_t u = - c\, q_x$ and, thus,
we look for solutions of the $1$D problem
\begin{equation}
\label{i12}
\left\{
\begin{aligned}
&   \mathrm{d}_x \left[ \partial\Phi (q_x)\right] + c\, q_x + f(q) = 0
    \quad\mbox{ for }\, x\in \mathbb{R} \,,
\\
& q_x\leq 0\,\mbox{ in }\, \mathbb{R} \quad\mbox{ and }\quad
  \lim_{x\to -\infty} q(x) = 1 \,,\quad
  \lim_{x\to \infty} q(x) = -1 \,,
\end{aligned}
\right.
\end{equation}
normalized by the condition
\begin{equation}
\label{i13}
  q(0) = \mu \,.
\end{equation}
As usual, we abbreviate the derivative
$q_x\equiv \mathrm{d}_x q\equiv q'$.
We show that the front speed $c$ as well as the solution $q$ of
\eqref{i12}, \eqref{i13} are unique, and
the asymptotic behavior of solutions to \eqref{i4}
is uniquely determined by \eqref{i10}, \eqref{i11}.
Since, by \eqref{i12}, we have
\begin{equation*}
    \int_{-\infty}^{\infty} |q_x| \,\mathrm{d}x
  = \int_{-\infty}^{\infty} (- q_x) \,\mathrm{d}x
  = - \lim_{x\to \infty} q(x) + \lim_{x\to -\infty} q(x)
  = 2
\end{equation*}
and the limits
$\lim_{x\to \pm\infty} \tilde\Phi(q_x)$ exist in $\mathbb{R}$ for
\begin{equation*}
  \tilde\Phi(z)\eqdef
  \int_0^z \partial^2\Phi(s)\, s \,\mathrm{d}s \,,\quad
    z\in \mathbb{R} \,,
\end{equation*}
where $\partial^2\Phi(s) > 0$ for $s\in \mathbb{R}\setminus \{ 0\}$,
so do the limits
$\lim_{x\to \pm\infty} q_x$.
Consequently, owing to $q_x\in L^1(\mathbb{R})$, we have also
\begin{equation}
\label{i12_x}
  q_x(x)\to 0 \quad\mbox{ as }\, x\to \pm\infty \,.
\end{equation}
Finally, multiplying eq.~\eqref{i12} by $q_x$ we get
\begin{equation*}
  \mathrm{d}_x \left[ \tilde\Phi (q_x)\right] + c\, |q_x|^2
  + \mathrm{d}_x F(q) = 0
    \quad\mbox{ for }\, x\in \mathbb{R} \,,
\end{equation*}
and integrating from $-\infty$ to $+\infty$ we arrive at
\begin{equation}
\label{e:c}
  c \int_{-\infty}^{\infty} |q_x|^2 \,\mathrm{d}x
  = \int_{-1}^{1} f(q) \,\mathrm{d}q = F(1) - F(-1) = F(1) \,.
\end{equation}

Note that condition~\eqref{int:f(u)}, i.e.,
$F(r) < F(1)$ for every $r\in (-1,1)$, forces $F(1)\geq F(-1) = 0$.
In order to exclude the {\em stationary solution\/}
$u(x,t) = q(x)$ with $c=0$, we will assume $F(1) > 0$,
in addition to \eqref{int:f(u)}, i.e.,
\begin{equation}
\label{int:F(1)>0}
  F(1) = \int_{-1}^{1} f(s) \,\mathrm{d}s > 0 \,.
\end{equation}

The paper is organized as follows.
In the next section (Section~\ref{s:prelim})
we recall some known results concerning solvability of
problem~\eqref{i1}, \eqref{i2} and state our main result in
Theorem~\ref{thm-sol}.
Section~\ref{s:trav_wave} is dedicated to the analysis of
the travelling wave problem \eqref{i12}, \eqref{i13}.
In particular, we show that problem \eqref{i12}, \eqref{i13}
admits a unique wave speed $c$ and a unique solution $q$
for a fairly general class of nonlinearities $\partial\Phi$ and $f$;
see Proposition~\ref{prop-TW}.
This means that no solution to this problem exists
for any other wave speed.
Finally, the convergence claimed in \eqref{i10}, \eqref{i11}
is established in Section~\ref{s:conv-TW}.

\section{Preliminaries, weak solutions, main result}
\label{s:prelim}

We shall say that a function
$u\colon \mathbb{R}^N\times \RR_+\to \RR$
is a \emph{weak solution\/} to problem \eqref{i1}, \eqref{i2} in
$\mathbb{R}^N\times \RR_+$
if it belongs to the class
\begin{align*}
& u\in L^\infty\left( \mathbb{R}^N\times (0, \infty) \right) ,\quad
  u\in C_{\mathrm{weak}}
       \left( [0,T]\to L^1_{\mathrm{loc}}(\mathbb{R}^N) \right)
    \quad\mbox{ for every }\, T > 0 \,,
\\
&   \mbox{ and }\quad
  \nabla_x u ,\; \partial\Phi(\nabla_x u)
  \in L^\infty\left( \mathbb{R}^N\times (0,T)\to \mathbb{R}^N\right) ,
\end{align*}
and satisfies the integral identity
\begin{equation}
\label{p1}
\begin{aligned}
&   \int_{\mathbb{R}^N} u(x,T)\, \psi(x,T) \,\mathrm{d}x
  - \int_{\mathbb{R}^N} u_0(x)\,    \psi(x,0)    \,\mathrm{d}x
\\
& = \int_0^{T} \int_{\mathbb{R}^N}
    \left[ u\, \partial_t\psi
        - \partial\Phi(\nabla_x u)\cdot \nabla_x\psi
        + f(u)\, \psi
    \right] \,\mathrm{d}x \,\mathrm{d}t
\end{aligned}
\end{equation}
required to hold for every test function
$\psi\in C_{\mathrm{c}}^{\infty}(\mathbb{R}^N\times \RR_+)$
and for every $T\geq 0$.
As usual,\hfil\break
$C_{\mathrm{c}}^{\infty}(\mathbb{R}^N\times \RR_+)$
denotes the space of all infinitely many times differentiable functions
$\psi\colon$\hfil\break
$\mathbb{R}^N\times \RR_+\to \mathbb{R}$
with compact support ($\subset \mathbb{R}^N\times \RR_+$).
For a weak solution to be well\--defined,
the function $u$ must obey the integrability conditions indicated above.

Analogously, we may define
the (weak) {\em sub-\/} and {\em super\-solutions\/} by changing
the equality sign ``$=$'' in eq.~\eqref{p1} to ``$\geq$'' and ``$\leq$'',
respectively, and taking there non\-negative test functions $\psi$ only.

In fact, with a help from
the ``regularity'' Proposition~\ref{prop-int_reg} below,
which guarantees
the (local H\"older-) continuity of a weak solution $u$ in
$\mathbb{R}^N\times (0,\infty)$,
we will construct
a \emph{viscosity solution\/} to problem \eqref{i1}, \eqref{i2} in
$\mathbb{R}^N\times \RR_+$.

\subsection{Existence of weak solutions, comparison principle, uniqueness}
${}$

$\bullet\;$
We assume that the function $\Phi\colon \mathbb{R}^N\to \mathbb{R}$
satisfies the following hypotheses:
\begin{align}
\label{p2}
\left\{
\begin{aligned}
  \Phi\colon \RR^N\to \RR \;\mbox{ is radially symmetric, i.e., }\;
  \Phi(Z)\equiv \varphi(|Z|) \;\mbox{ for every }\, Z\in \RR^N \,,
\\
  \Phi\in C^1(\RR^N)\cap C^2( \RR^N\setminus \{ 0\} ) \,,\quad
  \varphi\in C^1(\RR_+)\cap C^2((0,\infty)) \,,\quad
\end{aligned}
\right.
\\
\label{p3}
\left\{\quad
\begin{aligned}
& \varphi(0) = 0 \,,\quad \mathrm{d}_z\varphi (0) = 0 \,,
  \quad\mbox{ together with }\quad
\\
& \mathrm{d}_z\varphi(z) > 0
  \quad\mbox{ and }\quad
  \Lambda_1\leq
  \frac{ z\cdot \mathrm{d}^2_{z,z} \varphi(z) }{ \mathrm{d}_z \varphi(z) }
  \leq \Lambda_2
  \quad\mbox{ for all }\, z\in (0,\infty)\,,\quad
\\
& \mbox{ where }\, \Lambda_1, \Lambda_2 > 0
  \,\mbox{ are some positive constants. }
\end{aligned}
\right.
\end{align}

Moreover, there exists a continuous ``modulus of continuity'' function
$\omega\in C(\RR_+)$ with $\omega(0) = 0$, such that
\begin{equation}
\label{p3a}
\left\{\quad
\begin{aligned}
  \left| \mathrm{d}^2_{z,z} \varphi(s) - \mathrm{d}^2_{z,z} \varphi(y)
  \right|
  \leq \omega\genfrac{(}{)}{}0{|s-y|}{y}
       \cdot \mathrm{d}^2_{z,z} \varphi(y)
\\
  \quad\mbox{ for all $s,y\in (0,\infty)$ satisfying }\,
              |s-y| < \genfrac{}{}{}1{1}{2} y \,.
\end{aligned}
\right.
\end{equation}

$\bullet\;$
Besides the KPP condition \eqref{i3}, condition \eqref{int:f(u)}, and
$F(1) > 0$, we assume that
\begin{equation}
\label{p4}
  f\colon \RR\to \RR \quad\mbox{ is Lipschitz continuous. }
\end{equation}
%

\begin{remark}\label{rem-convex}\nopagebreak
\begingroup\rm
Condition \eqref{p3} implies that
$\varphi\colon \RR_+\to \RR$
is strictly monotone increasing and strictly convex.
Moreover, the Euclidean norm
$\vert\,\cdot\,\vert\colon \RR^N\to \RR$ is strictly convex
(even uniformly convex).
It is now an easy exercise to verify that also the function
$\Phi = \varphi\circ |\,\cdot\,|\colon \RR^N\to \RR$,
$\Phi(Z)\equiv \varphi(|Z|)$ for $Z\in \RR^N$,
must be {\em strictly convex\/}.
\endgroup
\end{remark}
\par\vskip 10pt

Hypotheses \eqref{p3}, \eqref{p3a} were introduced by
{\sc G.~M.\ Lieberman} \cite{Lieb3, Lieb4}.
In particular, both are satisfied for
a finite sum of $p$-Laplace operators with different exponents
$p = p_i\in (1,\infty)$.
Roughly speaking, hypothesis \eqref{p3} guarantees
{\it a~priori\/} bounds on
$\| \nabla_x u \|_{ L^{\infty}_{\mathrm{loc}} }$
in terms of $\| u \|_{ L^{\infty} }$ and,
in combination with \eqref{p3a},
also the H\"older continuity of $\nabla_x u$
for any bounded weak solution $u$ of \eqref{i1}, \eqref{i2}, see
{\sc G.~M.\ Lieberman} \cite{Lieb3, Lieb4}.
The implications of \eqref{p3}, \eqref{p3a}
on the structural properties of $\varphi$ as well as
other applications to the related elliptic problems are discussed by
{\sc Breit}, {\sc Stroffolini}, and {\sc Verde} \cite{BrStVe}.

\subsubsection{A priori bounds}

We quote the following crucial result on
\emph{interior} regularity estimates for equation \eqref{i1} from
{\sc G.~M.\ Lieberman} \cite{Lieb3}, \cite{Lieb4}:

\begin{proposition}\label{prop-int_reg}
Let\/ $f$ and\/ $\Phi$ satisfy hypotheses
\eqref{i3}, \eqref{int:f(u)}, and\/ \eqref{p2}--\eqref{p4}.
Assume that\/ $u$ is a weak solution of eq.~\eqref{i1} in
a bounded open space\--time cylinder
\begin{equation*}
  Q_{R,T} =
  \left\{ (x,t)\in \RR^N\times \RR_+\colon |x| < R \,,\  0 < t < T
  \right\} ,
\end{equation*}
for some $R,T\in (0,\infty)$,
belonging to the class
\begin{equation*}
  u\in L^{\infty}(Q_{R,T}) \,,\quad
  \nabla_x u\in L^{\infty}(Q_{R,T}\to \RR^N) \,.
\end{equation*}

{\bf (i)}
Then both $u, \nabla_x u$ are $\alpha$-H\"older continuous in the set
\begin{equation*}
  Q_{R,T,\delta} =
  \left\{ (x,t)\in Q_{R,T}\colon |x| < R - \delta \,,\ t > \delta
  \right\}
  \quad\mbox{ for any (sufficiently small) }\, \delta > 0 \,,
\end{equation*}
\begin{equation}
\label{est}
    \| u\|_{ C^{\alpha}(Q_{R,T,\delta}) }
  + \| \nabla_x u \|_{ C^{\alpha}(Q_{R,T,\delta}) }
  \leq M \,,
\end{equation}
where $\alpha$ and\/ $M$ depend solely on the set of parameters
\begin{equation*}
  \left[ R, T, \delta, \Lambda_1, \Lambda_2,
         \partial_z \varphi (1), \omega,
  \| u \|_{L^\infty(Q_{R,T})}, \| f(u) \|_{L^\infty(Q_{R,T})}
  \right] .
\end{equation*}

{\bf (ii)}
If, in addition,
\begin{equation*}
  u(\,\cdot\,, 0) = u_0
  \in C^{1+\beta}\left( \{ x\in \RR^N\colon |x|\leq R \} \right)
\end{equation*}
then \eqref{est} holds also for $\delta = 0$ in {\rm Part}~{\bf (i)},
i.e., in $Q_{R,T}$, with $\alpha$ and $M$ depending also on
$\beta$ and $\| u_0\|_{ C^{1+\beta} }$.
\end{proposition}

Proposition~\ref{prop-int_reg} yields an important corollary:

\begin{corollary}\label{cor-int_reg}
Let\/ $f$ and\/ $\Phi$ satisfy hypotheses
\eqref{i3}, \eqref{int:f(u)}, and\/ \eqref{p2}--\eqref{p4}.
Let\/ $u$ be a weak solution to problem \eqref{i1}, \eqref{i2}
in the set\/ $\RR^N\times (0,\infty)$ belonging to the class
\begin{equation*}
  u\in L^{\infty}( \RR^N\times (0,\infty) ) \,,\quad
  -1\leq u\leq 1 \,,\quad
  \nabla_x u\in L^{\infty}_{\mathrm{loc}}
    \left( \RR^N\times [0,\infty) \to \RR^N\right) ,
\end{equation*}
with the initial data
\begin{equation}
\label{p5}
\left\{
\begin{aligned}
& u_0\in BUC(\RR^N) \,,\ \nabla_x u_0\in BUC(\RR^N\to \RR^N) \,,\quad
  -1\leq u_0(x)\leq 1 \;\mbox{ for }\, x\in \RR^N \,,
\\
& | \nabla_x u_0(x) - \nabla_x u_0(y) |\leq L\, |x-y|^\beta
  \;\mbox{ for all }\, x,y\in \RR^N \;\mbox{ and certain }\,
  \beta\in (0,1] \,.
\end{aligned}
\right.
\end{equation}

Then there is a constant $M\in \RR_+$, depending solely on the quantities
\begin{equation*}
  \Lambda_1, \Lambda_2, \mathrm{d}_z\varphi(1), \omega, L, \beta,
  \| \nabla_x u_0 \|_{ L^\infty(\RR^N) } \,,
\end{equation*}
such that
\begin{equation}
\label{grad}
  | \nabla_x u(t,x) |\leq M
  \quad\mbox{ for all $x\in \RR^N$ and $t\in \RR_+$. }
\end{equation}
\end{corollary}

\subsubsection{Existence of solutions}

With the {\it a~priori} estimates stated in
Proposition~\ref{prop-int_reg} at hand, proving
\emph{existence\/} of a weak solution to
problem \eqref{i1}, \eqref{i2} is standard,
at least for smooth initial data.
More precisely, fixing the initial data $u_0$
in the regularity class \eqref{p5},
we may proceed in several steps, as follows:

{\it Step~$1$.}$\;$
Without any loss of generality, in \eqref{p3} we may assume
\begin{equation*}
  0 < \Lambda_1 < 1 < \Lambda_2 \,.
\end{equation*}
We perform a quadratic (Laplacian\--type) regularization of the function
$\Phi$ near the origin in $\RR^N$ by
replacing $\Phi$ by another $C^2$-smooth, radially symmetric function
$\Phi_{\alpha}\in C^2(\RR^N)$ ($\alpha > 0$), where
\begin{align}
\label{ap1}
& \quad
  \Phi_{\alpha}(Z)\equiv \varphi_{\alpha}(|Z|)
    \;\mbox{ for every }\, Z\in \RR^N \,,\quad
  \varphi_{\alpha}\in C^2(\RR_+) \,,\quad
\\
\label{ap2}
&
\left\{\quad
\begin{aligned}
& \varphi_{\alpha}(0) = 0 \,,\quad
  \mathrm{d}_z\varphi_{\alpha}(0) = 0 \,,\quad\mbox{ together with }\quad
\\
& \mathrm{d}_z\varphi_{\alpha}(z) > 0 \quad\mbox{ and }\quad
  \Lambda_1\leq
  \frac{ z\cdot \mathrm{d}^2_{z,z} \varphi_\alpha (z) }%
       { \mathrm{d}_z \varphi_\alpha (z) }
  \leq \Lambda_2 \quad\mbox{ for all }\, z > 0 \,,
\end{aligned}
\right.
\\
\label{ap2a}
&
\left\{\quad
\begin{aligned}
  \left| \mathrm{d}^2_{z,z} \varphi_{\alpha}(s)
       - \mathrm{d}^2_{z,z} \varphi_{\alpha}(y)
  \right|
  \leq \omega\genfrac{(}{)}{}0{|s-y|}{y}
       \cdot \mathrm{d}^2_{z,z} \varphi_{\alpha}(y)
\\
  \quad\mbox{ for all $s,y\in (0,\infty)$ satisfying }\,
              |s-y| < \genfrac{}{}{}1{1}{2} y \,,
\end{aligned}
\right.
\\
\label{ap3}
& \quad
  \frac{ z\cdot \mathrm{d}^2_{z,z} \varphi_{\alpha}(z) }%
       { \mathrm{d}_z \varphi_\alpha(z) }
  = 1 \quad\mbox{ for all }\, 0\leq z\leq \frac{1}{\alpha} \,,
    \qquad\mbox{ and }
\\
\label{ap4}
& \mathrm{d}_z \varphi_\alpha\to \mathrm{d}_z \varphi
  \quad\mbox{ as }\, \alpha\searrow 0 \,,
  \quad\mbox{ uniformly on compact subsets of }\, \RR_+ \,.
\end{align}
Here, the positive constants
$0 < \Lambda_1 < 1 < \Lambda_2$ and
the ``modulus of continuity'' function
$\omega\in C(\RR_+)$, with $\omega(0) = 0$,
are the same as in hypotheses \eqref{p3} and~\eqref{p3a};
in particular, all of them are independent from $\alpha > 0$.

We remark that hypothesis \eqref{ap3} is equivalent to
\begin{equation*}
  \mathrm{d}_z\, \log
  \genfrac{(}{)}{}0{ \mathrm{d}_z \varphi_\alpha(z) }{z}
  = 0 \quad\mbox{ for all }\, 0\leq z\leq \frac{1}{\alpha} \,,
\end{equation*}
which forces
$\varphi_\alpha(z) = \mathrm{const}_{\alpha}\cdot z^2$ for all
$0\leq z\leq \frac{1}{\alpha}$, with a positive constant.
Thus, $\varphi_\alpha$ is a qudratic regularization of the function
$\varphi\colon \RR_+\to \RR$ near zero.
This kind of regularization is typical in a construction of
a viscosity solution to a quasilinear parabolic problem.

{\it Step~$2$.}$\;$
Thanks to eq.~\eqref{ap3},
the resulting problem \eqref{i1}, \eqref{i2} with $\Phi$ replaced by
$\Phi_\alpha$ (defined in {\rm Step}~$1$)
is {\em uniformly\/} parabolic; thus,
by virtue of the standard theory for parabolic equations from
{\sc O.~A.\ Ladyzhenskaya}, {\sc V.~A.\ Solonnikov}, and
{\sc N.~N.\ Ural'tseva} \cite{LadySU},
it admits a unique (classical) solution $u_{\alpha}$
for any fixed $\alpha > 0$.

{\it Step~$3$.}$\;$
Since $f$ vanishes at $\pm 1$, the constant functions
$u(x,t)\equiv \pm 1$ are (classical) solutions of eq.~\eqref{i1}.
Thanks to
$-1\leq u_{\alpha}(\,\cdot\,,0) = u_0\leq 1$
in the initial condition \eqref{i2},
we may apply the classical version of
the (parabolic weak) comparison principle to deduce that also
\begin{equation}
\label{p6}
  -1\leq u_\alpha(t,x)\leq 1
    \quad\mbox{ for all }\, x\in \RR^N ,\ t\geq 0 .
\end{equation}

{\it Step~$4$.}$\;$
The family of classical solutions
$\{ u_{\alpha}\}_{\alpha > 0}$ satisfies the hypotheses of
Proposition~\ref{prop-int_reg} and Corollary~\ref{cor-int_reg},
with the parameters {\it independent\/} from~$\alpha$.
Consequently, it is easy to pass to the limit for $\alpha\searrow 0$,
at least for a suitable subsequence,
to deduce the following existence result:

\begin{proposition}\label{prop-exist}
Let the functions\/ $f$ and\/ $\Phi$ satisfy hypotheses
\eqref{i3}, \eqref{int:f(u)}, and\/ \eqref{p2}--\eqref{p4}.
Then for any initial data $u_0$ satisfying \eqref{p5},
the Cauchy problem \eqref{i1}, \eqref{i2} admits a weak solution $u$
in the class
\begin{equation*}
\begin{aligned}
& -1\leq u(t,x)\leq 1 \,,\quad
  |\nabla_x u(t,x)|\leq M
    \;\mbox{ for all }\, x\in \RR^N ,\ t\geq 0 ,
\\
& u,\ \nabla_x u \;\mbox{ belong to }\; C^{\kappa}(\mathcal{K}) ,\,
  \kappa\in (0,1) ,\;\mbox{ for any compact set }\,
  \mathcal{K}\subset \RR^N\times [0,\infty) \,.
\end{aligned}
\end{equation*}
\end{proposition}
\par\vskip 10pt

\subsubsection{Admissible weak solutions}

The construction procedure carried over in the preceding paragraph
inspires the following definition.

\begin{definition}\label{def-sol}\nopagebreak
\begingroup\rm
We say that $u$ is an {\em admissible weak solution\/}
to problem \eqref{i1}, \eqref{i2}, if there exists a sequence of
regularized functions $\varphi_{\alpha_n}$ ($n=1,2,3,\dots$)
enjoying properties \eqref{ap2}--\eqref{ap4} and a sequence of
initial data $u_{\alpha_n,0}$ belonging to
the regularity class \eqref{p5}, such that
\begin{align}
\label{dd1}
& u_{\alpha_n,0}\to u_0
  \quad\mbox{ uniformly on compact sets in }\, \RR^N \,,
  \;\mbox{ as }\, n\to \infty \,;
\\
\nonumber
& u_{\alpha_n}\to u
  \quad\mbox{ uniformly on compact sets in }\, \RR^N\times [0,\infty) \,,
\end{align}
where $\alpha_n\searrow 0$ as $n\nearrow \infty$, and
each $u_{\alpha_n}$ is the classical solution of
problem \eqref{i1}, \eqref{i2} corresponding to
$\Phi = \Phi_{\alpha_n}$ and $u_0 = u_{\alpha,0}$.
\endgroup
\end{definition}
\par\vskip 10pt

It can be shown that the admissible weak solutions coincide with
the standard {\em viscosity solutions\/} introduced for
continuous initial data by
{\sc Crandall}, {\sc Ishii}, and {\sc Lions} \cite{CIL}
as soon as $\varphi\in C^2[0,\infty)$.
As such, they satisfy the parabolic weak {\it comparison principle\/}
and, consequently, are uniquely determined by the initial data;
see {\sc Y.\ Giga} et al.~\cite{GiGoIsSa}.

Unfortunately, the singular case
$\mathrm{d}^2_{z,z} \varphi\to \infty$ for $z\searrow 0$,
that includes the $p$-Laplace operator with $1 < p < 2$,
does not fit into the framework of \cite{GiGoIsSa},
so that the mere definition of the concept of viscosity solution
requires some non\-trivial modifications, see
{\sc Juutinen}, {\sc Lindqvist}, and {\sc Manfredi} \cite{JuLiMa}.
Although the results and techniques used by
{\sc DiBenedetto} and {\sc Herrero} \cite{DibHer1, DibHer2}
provide the uniqueness of weak solutions for $p$-Laplace\--like operators
with $1 < p < 2$ and $p\geq 2$, respectively,
a general uniqueness theorem that would cover all cases allowed by
hypotheses \eqref{p3}, \eqref{p3a} does not seem to be easily available
in the existing literature.
Note, however, that a well\--posedness theory can be established by
the method of monotone operators
(see, e.g.,
 {\sc V.\ Barbu} \cite{Barbu} or {\sc H.\ Br\'ezis} \cite{Brezis})
as soon as the initial data approach one of the zeros of the function $f$
as $|x|\to \infty$.

\subsection{Main result}
\label{s:main}

We will show in Section~\ref{s:trav_wave}
that the traveling wave problem \eqref{i12}, \eqref{i13} admits
a {\em unique\/} solution pair $[q,c]$.
Accordingly, our main result may be stated as follows:

\begin{theorem}\label{thm-sol}
Let the functions\/ $f$ and\/ $\Phi$ satisfy hypotheses
\eqref{i3}, \eqref{int:f(u)}, $F(1) > 0$, and\/
\eqref{p2}--\eqref{p4}.
Assume that\/
$\{ u_{\eps}\}_{\eps > 0}$ is a family of admissible weak solutions of
the Cauchy problem \eqref{i4}, \eqref{i5},
with the initial data
$u_{\eps}(\,\cdot\,,0) = u_{\eps,0}\colon \mathbb{R}^N\to \RR$
satisfying
\begin{equation}
\label{p8}
\begin{aligned}
  -1\leq u_{\eps,0}(x)\leq 1 \quad\mbox{ and }\quad
  | \nabla_x u_{\eps,0}(x) |\leq \mathrm{const}_{\eps}
  \quad\mbox{ for all }\, x\in \RR^N \,,\quad\mbox{ and }
\\
  | \nabla_x u_{\eps,0}(x) - \nabla_x u_{\eps,0}(y) |
  \leq \mathrm{const}_{\eps}\cdot |x-y|^{\beta}
  \quad\mbox{ for all }\, x,y\in \RR^N \,,
\end{aligned}
\end{equation}
where $\beta\in (0,1)$ is a constant independent from $\eps > 0$.
Let $G_-$ and\/ $G_+$ be the sets introduced in eqs.\
\eqref{i7} and\/ \eqref{i8}, respectively, such that\/
eq.~\eqref{i9} holds, i.e.,
\begin{math}
  \overline{G}_- \cup \overline{G}_+ = \mathbb{R}^N \,.
\end{math}
Denote
\begin{math}
  \Gamma = \overline{G}_- \cap\overline{G}_+ \,.
\end{math}

Then, as\/ $\eps\searrow 0$, we have
\begin{align*}
  u_{\eps}\to -1
& \quad\mbox{ uniformly in compact subsets of }\,
  \left\{ (x,t)\colon \mathrm{dist}[x, \Gamma] > c t \right\} \,,
\\
  u_{\eps}\to 1
& \quad\mbox{ uniformly in compact subsets of }\,
  \left\{ (x,t)\colon \mathrm{dist}[x, \Gamma] < c t \right\} \,,
\end{align*}
with the \emph{front speed} $c$ being uniquely determined by eqs.\
\eqref{i12}, \eqref{i13}, $c > 0$,
where $\mathrm{dist}$ stands for the signed distance introduced
in~\eqref{DIST}.
\end{theorem}
\par\vskip 10pt

The remaining part of the paper is devoted to the proof of
Theorem~\ref{thm-sol}.

\section{Travelling waves}
\label{s:trav_wave}

In this section we establish the following result on
the existence and uniqueness of travelling waves in problem
\eqref{i12}, \eqref{i13}.

\begin{proposition}\label{prop-TW}
If the functions\/ $f$ and\/ $\Phi$ satisfy hypotheses
\eqref{i3}, \eqref{int:f(u)}, and\/ \eqref{p2}--\eqref{p4},
then problem \eqref{i12}, \eqref{i13} admits a unique solution $[q,c]$.
This solution satisfies also
\begin{align*}
& q_x\equiv \mathrm{d}_x q\in L^1(\RR)\cap L^{\infty}(\RR) \,,
\\
& q_x(x) < 0 \;\mbox{ for every $x\in \RR$ such that }\;
  -1 < q(x) < 1 \,,\quad\mbox{ and }
\\
& c = \left( \int_{-\infty}^{\infty} |q_x|^2 \,\mathrm{d}x\right)^{-1}
             \int_{-1}^{1} f(s) \,\mathrm{d}s \geq 0\,.
\end{align*}
In particular, we have $c > 0$ if and only if\/
$F(1) = \int_{-1}^{1} f(u) \,\mathrm{d}u > 0$.
\end{proposition}
\par\vskip 10pt

We will see that the {\it\bfseries proof\/} of this proposition
follows directly from a combination of
Lemma~\ref{lem-z_c(1):-infty} (existence) and
Proposition~\ref{prop-unique_c} (uniqueness) established below.
The continuous dependence of the travelling wave $q$ and the speed $c$
upon the given data,
combined with a standard compactness argument
(Arzel\`a\--Ascoli's theorem),
enables us to apply a continuity (convergence) result from
{\sc Ph.~Hartman}'s monograph
\cite[Theorem 2.1, p.~94]{Hartman}
to establish an approximation and continuity result
for problem \eqref{i12}, \eqref{i13} stated in the next lemma.

To formulate this result,
let us consider the following family of analogous problems
parametrized by $\alpha\in (0,1)$, for an unknown pair
$[q_{\alpha},c_{\alpha}]$
(cf.\ eqs.\ \eqref{ap1}--\eqref{ap4}):
\begin{equation}
\label{a:i12}
\left\{
\begin{aligned}
&   \mathrm{d}_x
  \left[ \partial\Phi_{\alpha} ( \mathrm{d}_x q_{\alpha} )\right]
  + c_{\alpha}\cdot \mathrm{d}_x q_{\alpha}
  + f_{\alpha}(q_{\alpha}) = 0
    \quad\mbox{ for }\, x\in \mathbb{R} \,,
\\
& \mathrm{d}_x q_{\alpha} \leq 0
    \,\mbox{ in }\, \mathbb{R} \quad\mbox{ and }\quad
  \lim_{x\to -\infty} q_{\alpha}(x) = b_{\alpha} \,,\quad
  \lim_{x\to \infty} q_{\alpha}(x) =  a_{\alpha} \,,
\end{aligned}
\right.
\end{equation}
where
$-\infty < a_{\alpha} < b_{\alpha} < \infty$,
with $q_{\alpha}$ normalized by the condition
\begin{equation}
\label{a:i13}
  q_{\alpha}(0) = \mu_{\alpha} \,,\quad
    a_{\alpha} < \mu_{\alpha} < b_{\alpha} \,.
\end{equation}
The existence and uniqueness of the pair
$[q_\alpha, c_\alpha]$ follow in the same way as those of $[q,c]$, cf.\
{\rm Proposition~\ref{prop-TW}} above.

\begin{lemma}\label{lem-TW}
Let\/
$\{ f_{\alpha}\}_{\alpha > 0}$
be a family of uniformly Lipschitz\--continuous functions
$f_{\alpha}\colon \RR\to \RR$, such that\/
\begin{align}
\label{a:Lip}
& L\eqdef
  \sup_{\alpha > 0} \| f_{\alpha}'\|_{ L^{\infty}(\RR) }
  < \infty \,;
\\
\label{a:sign_f}
&
\left\{
\begin{aligned}
& -\infty < a_{\alpha} < \mu_{\alpha} < b_{\alpha} < \infty \,,\quad
  f_{\alpha} (a_{\alpha}) = f(\mu_{\alpha}) = f(b_{\alpha}) = 0
\\
&   \quad\mbox{ with }\quad
  f_{\alpha} < 0 \,\mbox{ in }\, (a_{\alpha}, \mu_{\alpha}) \,,\quad
  f_{\alpha} > 0 \,\mbox{ in }\, (\mu_{\alpha}, b_{\alpha}) \,;
\end{aligned}
\right.
\\
\label{a:int_f}
&
\left\{
\begin{aligned}
& F_{\alpha}(r)\eqdef
  \int_{ a_{\alpha} }^r f_{\alpha}(s) \,\mathrm{d}s
    \quad\mbox{ satisfies }
\\
& F_{\alpha}( b_{\alpha} ) - F_{\alpha}(r) =
  \int_r^{ b_{\alpha} } f_{\alpha}(s) \,\mathrm{d}s > 0
    \quad\mbox{ whenever }\, a_{\alpha} < r < b_{\alpha} \,;
\end{aligned}
\right.
\\
\label{a:conv_f}
&
\left\{
\begin{aligned}
& a_{\alpha}\to -1 \,,\quad \mu_{\alpha}\to \mu \,,\quad
  b_{\alpha}\to 1 \,,\quad\mbox{ and }\quad
\\
& f_{\alpha}\to f \;\mbox{ locally uniformly in }\, C(\RR)
    \quad\mbox{ as }\, \alpha\searrow 0 \,,
\end{aligned}
\right.
\end{align}
and the limit function $f$ satisfies all conditions in \eqref{i3},
$-1 < \mu < 1$, together with condition \eqref{int:f(u)}.
Similarly, we assume that the family
$\{ \Phi_{\alpha}\}_{\alpha > 0}$
satisfies all hypotheses in \eqref{ap1}--\eqref{ap4}.
Finally, let\/
$[q_{\alpha}, c_{\alpha}]$ $(\alpha > 0)$
be the uniquely determined family of solutions
to problem \eqref{a:i12}, \eqref{a:i13} satisfying
\begin{equation*}
  \lim_{x\to -\infty} q_{\alpha}(x) = b_{\alpha} \,,\quad
  \lim_{x\to \infty}  q_{\alpha}(x) = a_{\alpha} \,,
    \quad\mbox{ and }\quad
  q_{\alpha} (0) = \mu_{\alpha} \,,
\end{equation*}
by\/ {\rm Proposition~\ref{prop-TW}}.
Then we have
\begin{align*}
&   \| \mathrm{d}_x q_{\alpha} \|_{ L^1(\RR) }
  + \| \mathrm{d}_x q_{\alpha} \|_{ L^{\infty}(\RR) }
  \leq \mathrm{const} \quad\mbox{ together with }
\\
& q_{\alpha}\to q \;\mbox{ uniformly in }\, C(\RR)
    \quad\mbox{ and }\quad
  c_{\alpha}\to c \quad\mbox{ as }\, \alpha\searrow 0 \,,
\end{align*}
where $[q,c]$ is the unique solution of \eqref{i12}, \eqref{i13}
corresponding to the limit functions $f$ and\/ $\Phi$,
by\/ {\rm Proposition~\ref{prop-TW}}.
\end{lemma}
\par\vskip 10pt

We begin with the proof of Proposition~\ref{prop-TW}.
Let us recall that we investigate monotone (decreasing) travelling waves
in the degenerate second\--order para\-bolic problem \eqref{i1}
of a ``generalized'' Fisher\--KPP type.
This task reduces to finding travelling waves in
the following degenerate second\--order para\-bolic problem
reduced to one space dimension:
\begin{equation}
\label{e:FKPP}
\left\{
\begin{aligned}
    \partial_t u
& = \partial_x
    \left( \partial\Phi (\partial_x u) \right)
  + f(u) \,,\quad (x,t)\in \RR^1\times \RR_+ \,,
\\
  u(x,t)&= q(x-ct) \quad\mbox{ for some constant }\, c\in \RR \,.
\end{aligned}
\right.
\end{equation}
Here, $\Phi\colon \RR^1\to \RR$
is an continuously differentiable, even convex function satisfying
all hypotheses in \eqref{p2} and \eqref{p3}, with the derivative
$\partial\Phi\equiv \Phi'\colon \RR^1\to \RR$.

For example, we may take
$\Phi(z) = \frac{1}{p} |z|^p$ for $z\in \RR$, where
$1 < p < \infty$ is a fixed number, so that
$\partial\Phi(z)\equiv \Phi'(z) = |z|^{p-2} z$ and
$\partial^2\Phi(z)\equiv \Phi''(z) = (p-1) |z|^{p-2}$ for $z\in \RR$.

Recall that $f\colon \RR\to \RR$ is assumed to be Lipschitz\--continuous,
by \eqref{p4}, and, most importantly, it satisfies
the KPP condition \eqref{i3}, that is,
$f(\pm 1) = f(\mu) = 0$ for some $-1 < \mu < 1$, together with
$f(s) < 0$ for every $s\in (-1,\mu)$,
$f(s) > 0$ for every $s\in (\mu,1)$, and also
condition \eqref{int:f(u)}, specifically
\begin{equation*}
\nonumber
  F(1) - F(r) = \int_r^1 f(s) \,\mathrm{d}s > 0
    \quad\mbox{ whenever }\, -1 < r < 1 \,.
\end{equation*}
Furthermore, by Proposition~\ref{prop-TW},
we have $F(1) > 0$ if and only if $c > 0$.

\begin{remark}\label{rem-double-well}\nopagebreak
\begingroup\rm
An important special case of the reaction function $f$ is
$f(s) = F'(s)$ where
$-F\colon \RR\to \RR$ is a ``generalized'' {\em double\--well potential\/}
(\cite[eq.\ (0.2)]{BBS}),
such that
\begin{equation*}
\left\{
\begin{aligned}
  F'(s)\equiv f(s)
& = 2 (s - \mu) (1 - s^2) = - 2 (s+1) (s - \mu) (s-1)
  \quad\mbox{ for }\, s\in \RR \,;
\\
  F(-1) &= 0 \,,
\end{aligned}
\right.
\end{equation*}
whence
\begin{equation}
\label{e:F}
\left\{
\begin{aligned}
  - F(s)
  = (s^2 - 1) (s - \mu)^2
  - \frac{1}{2}\, (s - \mu)^4
  - \frac{2}{3}\,\mu (s - \mu)^3 - F(\mu)
\\
    \quad\mbox{ for }\, s\in \RR \,,\quad\mbox{ where }\quad
  F(\mu) = {}
  - \frac{1}{2}\, (1 + \mu)^3 \left( 1 - \frac{\mu}{3}\right) < 0 \,.
\end{aligned}
\right.
\end{equation}
We have also
$F(1) = - F(-\mu) + F(\mu) = {}- \frac{8}{3}\, \mu$
or
\begin{equation*}
\begin{aligned}
  F(1)
& = \int_{-1}^{1} f(s) \,\mathrm{d}s
  = 2\int_{0}^{1} [ (s - \mu) + (- s - \mu) ] (1 - s^2) \,\mathrm{d}s
\\
& = {}- 4\mu\int_{0}^{1} (1 - s^2) \,\mathrm{d}s
  = {}- \frac{8}{3}\, \mu \,.
\end{aligned}
\end{equation*}
In particular, condition~\eqref{int:f(u)} holds if and only if
$\mu\leq 0$.
\endgroup
\end{remark}
\par\vskip 10pt

Assuming that the travelling wave takes the form
$u(x,t) = q(x-ct)$, $(x,t)\in \RR\times \RR_+$, with
$q\colon \RR\to \RR$ being continuously differentiable and satisfying
$q'(x) < 0$ at every point $x\in \RR$ such that $-1 < q(x) < 1$, below,
we are able to find a {\it\bfseries first integral\/} for
the second\--order equation for $q$; cf.\
eq.~\eqref{i12}:
\begin{equation}
\label{eq:FKPP}
    \mathrm{d}_x
    \left( \partial\Phi (\mathrm{d}_x q) \right)
  + c\cdot \mathrm{d}_x q + f(q)
  = 0 \,,\quad x\in \RR \,.
\end{equation}

We will use the following (abuse of) notation exclusively
throughout the remaining part of this section;
it will not intervene with the notation for the function
$\varphi = \varphi(|Z|) = \Phi(Z)$ of $Z\in \mathbb{R}^N$
introduced in Section~\ref{s:prelim}, eq.~\eqref{p2}.
This time, let us denote
$\varphi = \partial\Phi\equiv \Phi'$.
By the properties \eqref{p2} and \eqref{p3} of $\Phi$ recalled above,
$\varphi\colon \RR\to \RR$ is a continuous,
strictly monotone increasing, odd function; hence, $\varphi(0) = 0$.
Moreover, we have
$\varphi(s)\to \pm\infty$ as $s\to \pm\infty$, respectively.
Consequently,
\begin{equation}
\label{e:Phi'=phi}
  \Phi(s)\eqdef \int_0^s \varphi(\xi) \,\mathrm{d}\xi \,,
  \quad\mbox{ for }\, s\in \RR \,,
\end{equation}
is a continuously differentiable, strictly convex, even function, with
$\Phi(0) = 0$ and
$\Phi(s) / |s|$ $\to +\infty$ for $s\to \pm\infty$.
We denote by
$\Psi\colon \RR\to \RR$
the {\em convex conjugate\/} function associated with $\Phi$, that is,
\begin{equation}
\label{e:Psi/Phi}
  \Psi(t)\eqdef \sup_{s\in \RR} \left( st - \Phi(s)\right)
  \quad\mbox{ for }\, t\in \RR \,.
\end{equation}
Then, by the general theory for pairs of convex conjugate functions, also
$\Psi$ is continuously differentiable, strictly convex, and even, with
$\Psi(0) = 0$ and
$\Psi(t) / |t|$ $\to +\infty$ for $t\to \pm\infty$.
Its derivative
$\psi\eqdef \Psi'\colon \RR\to \RR$
is continuous, strictly monotone increasing, and odd.
Moreover, the functions
$\varphi, \psi\colon \RR\to \RR$ are each other's inverse, i.e.,
$\psi = \varphi_{-1}$ and $\varphi = \psi_{-1}$.
We refer the reader to the monograph by
{\sc I.\ Ekeland} and {\sc R.\ Temam} \cite[Part~1]{Eke-Temam}
for details about convex conjugate functions.

Following the main ideas from
{\sc R.\ Engui\c{c}a}, {\sc A.\ Gavioli}, and {\sc L.\ Sanchez}
\cite[Sect.~1]{EGSanchez},
we make the substitution
\begin{equation*}
  v\eqdef {}- \varphi(\mathrm{d}_x q) > 0 \,,\quad\mbox{ where }\quad
  \mathrm{d}_x q\equiv q_x\equiv \frac{\mathrm{d}q}{\mathrm{d}x} < 0 \,,
\end{equation*}
whence
\begin{equation}
\label{e:dU/dx}
    \mathrm{d}_x q = {}- \varphi_{-1}(v) = {}- \psi(v) < 0 \,,
\end{equation}
and consequently look for $v = v(q)$ as a function of
$q\in (-1,1)$ that satisfies the following differential equation
obtained from eq.~\eqref{eq:FKPP}:
\begin{equation*}
{}- \frac{\mathrm{d}v}{\mathrm{d}q}\cdot \frac{\mathrm{d}q}{\mathrm{d}x}
  + c\, \frac{\mathrm{d}q}{\mathrm{d}x} + f(q)
  = 0 \,,
    \quad x\in \RR \,,
\end{equation*}
that is, with a help from \eqref{e:dU/dx},
\begin{equation}
\label{eq:FKPP:V(U)}
    \frac{\mathrm{d}v}{\mathrm{d}q}\cdot \psi(v)
  - c\, \psi(v) + f(q) = 0 \,,
    \quad q\in (-1,1) \,.
\end{equation}
Finally, we make the substitution $y\eqdef \Psi(v) > 0$,
thus arriving at
\begin{equation*}
    \frac{\mathrm{d}y}{\mathrm{d}q}
  - c\, \psi\left( \Psi_{-1}(y) \right) + f(q) = 0 \,,
    \quad q\in (-1,1) \,.
\end{equation*}
Here,
\begin{math}
  \Psi_{-1}\equiv \left( \Psi\vert_{\RR_+} \right)_{-1}
  \colon \RR_+\to \RR
\end{math}
stands for the inverse function of
$\Psi$ restricted to the domain $\RR_+\eqdef [0,\infty)$ and, thus,
denoted by $\Psi\vert_{\RR_+}$.
In order to avoid possible confusion between
the {\it\bfseries unknown function\/} $q(x)$ of $x\in \RR$ and
the {\it\bfseries independent variable\/} $q\in (-1,1)$,
we prefer to replace the latter by $r\in (-1,1)$.
This means that the unknown function
$y\colon (-1,1)\to (0,\infty)$ of $r$,
\begin{equation}
\label{e:y=V^p'}
  y = \Psi(v) = \Psi\left( \varphi(|\mathrm{d}_x q|) \right) > 0 \,,
\end{equation}
must satisfy the following differential equation:
\begin{equation}
\label{eq:FKPP:y(r)}
    \frac{\mathrm{d}y}{\mathrm{d}r}
  - c\, \psi\left( \Psi_{-1}(y^{+}) \right) + f(r) = 0 \,,
    \quad r\in (-1,1) \,,
\end{equation}
where $y^{+}\eqdef \max\{ y,\, 0\}$ for $y\in \RR$.
Since we require that $q = q(x)$ be sufficiently smooth,
at least continuously differentiable, with
$q_x(x)\equiv q'(x)\to 0$ as $x\to \pm\infty$,
the function $y = y(r)$ must satisfy the boundary conditions
\begin{equation}
\label{bc:FKPP:y(r)}
  y(-1) = y(1) = 0 \,.
\end{equation}
Recalling the substitution $y = \Psi(v) > 0$ for $v > 0$, i.e.,
$v = \Psi_{-1}(y)$, from eq.~\eqref{eq:FKPP:V(U)}
we deduce the following equivalent form of eq.~\eqref{eq:FKPP:y(r)}
for the unknown function $v = v(r)$,
\begin{equation}
\label{eq:FKPP:v(r)}
    \frac{\mathrm{d}v}{\mathrm{d}r}
  - c + \frac{f(r)}{\psi(v)} = 0 \,,
    \quad r\in (-1,1) \,.
\end{equation}
Furthermore, boundary conditions \eqref{bc:FKPP:y(r)}
for $y = y(r)$ become
\begin{equation}
\label{bc:FKPP:v(r)}
  v(-1) = v(1) = 0 \,.
\end{equation}

The following remark on the value of $F(1)$ (${}\geq 0$) is in order.
We recall that $f(r)$ satisfies the KPP condition \eqref{i3},
together with condition \eqref{int:f(u)}.

\begin{remark}\label{rem-sol:F(r)}\nopagebreak
\begingroup\rm
Since the integrand $f\colon (-1,1)\to \RR$ in the function $F(r)$,
used in~\eqref{int:f(u)} for $r\in (-1,1)$,
is continuous and absolutely integrable over $(-1,1)$,
we conclude that
$F\colon [-1,1]\to \RR$ is absolutely continuous.
In particular, ineq.~\eqref{int:f(u)} forces $F(1)\geq 0$.
We will see later that the case $F(1) = 0$ guarantees
the existence of a {\em stationary solution\/}
to problem~\eqref{e:FKPP}, i.e., $c=0$, whereas
the case $F(1) > 0$ renders a {\em travelling wave\/}, i.e., $c\neq 0$;
more precisely, $c>0$, cf.\
Proposition~\ref{prop-TW} above.
Indeed, both,
the stationary solution (for $c = 0$) and
the travelling wave (for $c > 0$)
will be obtained from eq.~\eqref{eq:FKPP:y(r)}
by means of the transformation defined by eqs.\
\eqref{e:dU/dx} and \eqref{e:y=V^p'}.

Recalling Remark~\ref{rem-double-well},
for the quartic double\--well potential $-F$ given by eq.~\eqref{e:F}
we have
\begin{equation*}
  c = \left( \int_{-\infty}^{\infty} |q_x|^2 \,\mathrm{d}x\right)^{-1}
             \int_{-1}^{1} f(s) \,\mathrm{d}s
    = {}- \frac{8}{3}\, \mu
      \left( \int_{-\infty}^{\infty} |q_x|^2 \,\mathrm{d}x\right)^{-1} \,.
\end{equation*}
\endgroup
\end{remark}
\par\vskip 10pt

In order to investigate equation
\eqref{eq:FKPP:y(r)} (and \eqref{eq:FKPP:v(r)}, as well),
we begin with the following more general differential equation
than eq.~\eqref{eq:FKPP:y(r)}, namely,
\begin{equation}
\label{gen:FKPP:y(r)}
    \frac{\mathrm{d}y}{\mathrm{d}r} - c\, H(y^{+}) = - f(r) \,,
    \quad r\in (-1,1) \,,
\end{equation}
where $H\colon \RR_+\to \RR$ is a continuous,
strictly monotone increasing function with $H(0) = 0$.
In eq.~\eqref{eq:FKPP:y(r)}, this means
$H = \psi\circ \Psi_{-1}$ on $\RR_+$.

The positivity condition \eqref{int:f(u)} for $r\in (-1,1)$
starts from the terminal value of $r$ ($r=1$).
This would cause serious difficulties with notation and
the uniqueness for a number of initial value problems
that we are going to treat; namely,
we would be forced to treat them as terminal value problems.
Therefore, we make the substitution
$z(r) = y(-r)$ for $-1\leq r\leq 1$ (reflection about $0$)
for the unknown function $y$ and, consequently,
look for a continuously differentiable solution
$z\colon [-1,1]\to \mathbb{R}$
to the following Dirichlet boundary value problem equivalent to
eq.~\eqref{gen:FKPP:y(r)},
\begin{equation}
\label{gen:FKPP:z(r)}
    \frac{\mathrm{d}z}{\mathrm{d}r} + c\, H(z^{+}) = f(-r) \,,
    \quad r\in (-1,1) \,;\qquad
  z(-1) = z(1) = 0 \,,
\end{equation}
where $c\in \RR$ is also an unknown constant to be determined.
As usual, we apply the shooting method for solutions
$z\colon [-1,1]\to \mathbb{R}$ using the initial condition $z(-1) = 0$.
We determine the constant $c\in \RR$ such that also $z(1) = 0$ hold.

The following comparison lemma is standard; cf.\
{\sc Ph.~Hartman} \cite[Corollary 4.2, p.~27]{Hartman}.

\begin{lemma}\label{lem-sub-sup}
Let\/ $c\in \RR$ satisfy\/ $c\geq 0$ and assume that\/
$z_1, z_2\colon [a,b]\to \RR$ are two absolutely continuous functions
on some interval\/
$[a,b]\subset \RR$, such that\/ $-1\leq a < b\leq 1$,
and the following inequality holds for almost every $r\in (a,b)$:
\begin{equation}
\label{gen:FKPP:z_1<z_2}
    \frac{\mathrm{d}z_1}{\mathrm{d}r} + c\, H(z_1^{+})
  \leq
    \frac{\mathrm{d}z_2}{\mathrm{d}r} + c\, H(z_2^{+}) \,.
\end{equation}
If\/ $z_1(a)\leq z_2(a)$ then
$z_1(r)\leq z_2(r)$ holds for all\/ $r\in [a,b]$.
\end{lemma}

\par\vskip 10pt
\proof
Rewriting ineq.~\eqref{gen:FKPP:z_1<z_2}
for the difference $z_1 - z_2$, we have
\begin{equation*}
    \frac{\mathrm{d}}{\mathrm{d}r} (z_1 - z_2)
  + c\left( H(z_1^{+}) - H(z_2^{+}) \right) \leq 0
    \quad\mbox{ for a.e. }\, r\in (a,b) \,.
\end{equation*}
Multiplying this difference by $(z_1 - z_2)^{+}$, we arrive at
\begin{equation}
\label{gen:FKPP:z_1-z_2}
\begin{aligned}
  \frac{1}{2}\cdot
    \frac{\mathrm{d}}{\mathrm{d}r} \left[ (z_1 - z_2)^{+} \right]^2
  + c\left( H(z_1^{+}) - H(z_2^{+}) \right) (z_1 - z_2)^{+}
  \leq 0
\\
    \quad\mbox{ for a.e. }\, r\in (a,b) \,.
\end{aligned}
\end{equation}
Taking into account that $c\geq 0$ and
\begin{math}
  s\mapsto H(s^{+})\colon \RR\to \RR
\end{math}
is a monotone increasing function, with $H(s^{+}) = 0$ for $s\leq 0$,
we conclude that the second summand on the left\--hand side of
eq.~\eqref{gen:FKPP:z_1-z_2}
is nonnegative, which yields
\begin{equation*}
  \frac{1}{2}\cdot
    \frac{\mathrm{d}}{\mathrm{d}r} \left[ (z_1 - z_2)^{+} \right]^2
  \leq 0
    \quad\mbox{ for a.e. }\, r\in (a,b) \,.
\end{equation*}
Consequently,
\begin{math}
  r\mapsto \left[ (z_1(r) - z_2(r))^{+} \right]^2
  \colon [a,b]\to \RR
\end{math}
is a nonnegative, monotone decreasing function
that vanishes at $r=a$; hence,
it must vanish identically on the whole of $[a,b]$.

The lemma is proved.
\qed
\par\vskip 10pt

Lemma~\ref{lem-sub-sup} has the following easy, but very useful corollary.

\begin{corollary}\label{cor-sub-sup}
Let\/ $c_1, c_2\in \RR$ satisfy\/ $c_1\geq c_2$ and\/ $c_1\geq 0$.
Assume that\/
$z_1, z_2\colon [a,b]\to \RR$ are two absolutely continuous functions
on some interval\/
$[a,b]\subset \RR$, such that\/ $-1\leq a < b\leq 1$,
and the following inequalities hold for almost every $r\in (a,b)$:
\begin{align}
\label{gen:FKPP:z_1}
    \frac{\mathrm{d}z_1}{\mathrm{d}r} + c_1\, H(z_1^{+})\leq f(-r) \,,
\\
\label{gen:FKPP:z_2}
    \frac{\mathrm{d}z_2}{\mathrm{d}r} + c_2\, H(z_2^{+})\geq f(-r) \,.
\end{align}
If\/ $z_1(a)\leq z_2(a)$ then
$z_1(r)\leq z_2(r)$ holds for all\/ $r\in [a,b]$.
\end{corollary}

\par\vskip 10pt
\proof
Inequalities \eqref{gen:FKPP:z_1} and \eqref{gen:FKPP:z_2}
guarantee that \eqref{gen:FKPP:z_1<z_2}
holds with any $c\in \RR$ such that $c\geq 0$ and
$c_2\leq c\leq c_1$, e.g., with $c = c_1$.
Lemma~\ref{lem-sub-sup} yields the desired inequality.
\qed
\par\vskip 10pt

Corollary~\ref{cor-sub-sup} shows that the initial value problem for
eq.~\eqref{gen:FKPP:z(r)} with the initial condition
$z(-1) = 0$ possesses a {\em unique\/} (absolutely continuous) solution
$z\equiv z_c\colon [-1,1]\to \RR$,
whenever $c\geq 0$ is a fixed number.
A standard combination of compactness
(Arzel\`a\--Ascoli's theorem)
and uniqueness from
{\sc Ph.~Hartman} \cite[Theorem 2.1, p.~94]{Hartman}
guarantees that the solution mapping
$c\mapsto z_c\colon \RR_+\to C([-1,1])$
is continuous.

Note that, for any $c\leq 0$,
an arbitrary (possibly nonunique) solution
$z\colon [-1,1]\to \RR$ to eq.~\eqref{gen:FKPP:z(r)}
with $z(-1) = 0$ must satisfy
\begin{equation*}
    z(r)\geq \int_{-1}^{r} f(-s) \,\mathrm{d}s
  = \int_{-r}^{1} f(s') \,\mathrm{d}s' = F(1) - F(-r) > 0
    \quad\mbox{ for every }\, r\in (-1,1) \,,
\end{equation*}
by a simple integration of this equation over the interval $[-1,1]$
followed by ineq.~\eqref{int:f(u)}.
Assuming $F(1) > 0$, we get also
\begin{math}
  z(1)\geq \int_{-1}^{1} f(-s) \,\mathrm{d}s = F(1) > 0 .
\end{math}

Consequently, assuming ineq.~\eqref{int:f(u)} and $F(1) > 0$,
let us define
\begin{equation}
\label{def:c^*}
  c^{\ast}\eqdef \sup\left\{ c\in \RR_+\colon z_c(1) > 0\right\} \,.
\end{equation}
Clearly, $c^{\ast}\geq 0$.
As expected, we will show that precisely $c^{\ast}$
is the desired critical value of the constant~$c$, and $c^{\ast} > 0$.
From the continuity of the mapping
$c\mapsto z_c(1)\colon \RR_+\to \RR$
combined with $z_0(1) > 0$ (for $c=0$), we deduce that
either $c^{\ast} = +\infty$, or else
$0 < c^{\ast} < +\infty$ in which case $z_{c^{\ast}}(1) = 0$
and, consequently, $c^{\ast}$ is the desired critical value.
In what follows we exclude the former case, $c^{\ast} = +\infty$,
which would force $z_c(1) > 0$ for every $c\in \RR_+$.

Lemma~\ref{lem-sub-sup} has another important corollary
which, under stronger hypotheses on $H$, streng\-thens the conclusion of
Corollary~\ref{cor-sub-sup}.

\begin{corollary}\label{cor-sub<sup}
Let\/ $c_1, c_2\in \RR$ satisfy\/ $c_1 > c_2$ and\/ $c_1\geq 0$.
Assume that\/
$z_1, z_2\colon [a,b]\to \RR$ are two absolutely continuous functions
on some interval\/
$[a,b]\subset \RR$, such that\/ $-1\leq a < b\leq 1$,
and the following equations hold for almost every $r\in (a,b)$:
\begin{align}
\label{eq:FKPP:z_1}
    \frac{\mathrm{d}z_1}{\mathrm{d}r} + c_1\, H(z_1^{+}) = f(-r) \,,
\\
\label{eq:FKPP:z_2}
    \frac{\mathrm{d}z_2}{\mathrm{d}r} + c_2\, H(z_2^{+}) = f(-r) \,.
\end{align}
\begin{itemize}
\item[{\rm (i)}]
If\/ $z_1(a)\leq z_2(a)$
then precisely one of the following two alternatives occurs:
either\/
$z_1(r) < z_2(r)$ or else\/
$z_1(r) = z_2(r)\leq 0$ for every\/ $r\in (a,b)$.
\item[{\rm (ii)}]
If\/ $-\mu\leq a < b\leq 1$ and $z_1(a)\leq z_2(a)\leq 0$, then
$z_2(r) - z_1(r) = z_2(a) - z_1(a)$
holds for all\/ $r\in [a,b]$.
More precisely, we have
\begin{equation}
\label{eq:z_1,2(r)<0}
\begin{aligned}
  z_i(r) = z_i(a) + \int_a^r f(-s) \,\mathrm{d}s
         = z_i(a) + \int_{-r}^{-a} f(s) \,\mathrm{d}s
         < z_i(a)\leq 0
\\
    \quad\mbox{ for all }\, r\in (a,b] \,;\quad i=1,2 \,.
\end{aligned}
\end{equation}
\item[{\rm (iii)}]
Assume that\/
$H = \psi\circ \Psi_{-1}$ on $\RR_+$.
If\/ $-\mu\leq a < b\leq 1$, $z_1(a)\leq z_2(a)$, and\/
$z_1(r) > 0$ for all\/ $r\in (a,b)$, then
$z_1(r) < z_2(r)$ holds for all\/ $r\in (a,b]$.
\end{itemize}
\end{corollary}

\par\vskip 10pt
\proof
Part~{\rm (i)}:$\;$
It follows from
eqs.~\eqref{eq:FKPP:z_1} and~\eqref{eq:FKPP:z_2}
that both functions
$z_1, z_2\colon [a,b]\to \RR$
are continuously differentiable.
Assuming $z_1(a)\leq z_2(a)$,
from Corollary~\ref{cor-sub-sup} we deduce
$z_1(r)\leq z_2(r)$ for all\/ $r\in [a,b]$.
Now suppose there is some $r_0\in (a,b)$ such that
$z_1(r_0) = z_2(r_0) > 0$.
Consequently, we have also
$\frac{\mathrm{d}}{\mathrm{d}r} (z_2 - z_1)(r_0) = 0$.
We insert the equalities
$z_1(r_0) = z_2(r_0) > 0$ and $z_1'(r_0) = z_2'(r_0)$
into eqs.~\eqref{eq:FKPP:z_1} and~\eqref{eq:FKPP:z_2},
thus arriving at
$c_1\, H(z_1(r_0)) = c_2\, H(z_2(r_0))$.
Since
$H(z_1(r_0)) = H(z_2(r_0)) > 0$, we conclude that $c_1 = c_2$
which contradicts our hypothesis $c_1 > c_2$.

Part~{\rm (ii)}:$\;$
This part is derived directly from eq.~\eqref{gen:FKPP:z(r)}
considered for $r\in (a,b)$ with an unknown function
$z\colon [a,b]\to \RR$ that is assumed to be absolutely continuous.
Assuming
$-\mu\leq a < b\leq 1$ and $z(a)\leq 0$,
for every $r\in (a,b)$ we have $-1\leq -b < -r < -a\leq \mu$
which yields $f(-r) < 0$ and, consequently,
by eq.~\eqref{gen:FKPP:z(r)},
$z\colon [a,b]\to \RR$
is a strictly monotone decreasing function satisfying
\begin{equation}
\label{eq:z(r)<0}
\begin{aligned}
  z(r) = z(a) + \int_a^r f(-s) \,\mathrm{d}s
       = z(a) + \int_{-r}^{-a} f(s) \,\mathrm{d}s
         < z(a)\leq 0
\\
    \quad\mbox{ for all }\, r\in (a,b] \,.
\end{aligned}
\end{equation}
Notice that, in this part, the value of the constant $c\in \RR$
is completely irrelevant.

Part~{\rm (iii)}:$\;$
Next, assume
$-\mu\leq a < b\leq 1$, $z_1(a)\leq z_2(a)$, and
$z_1(r) > 0$ for all $r\in (a,b)$.
Then also $z_2(r)\geq z_1(r) > 0$ holds for all $r\in (a,b)$,
by Corollary~\ref{cor-sub-sup}.
Therefore, we can take advantage of the substitution
$y = \Psi(v)$ defined in eq.~\eqref{e:y=V^p'},
combined with $z(r) = y(-r)$ and $w(r) = v(-r)$ for $r\in [-1,1]$,
and use eq.~\eqref{eq:FKPP:v(r)} for $v(r)$ in place of
\eqref{eq:FKPP:y(r)} for $y(r)$.
For the unknown function $w = w(r) = v(-r)$ in place of $v$,
$w = \Psi_{-1}(z)$, eq.~\eqref{eq:FKPP:v(r)} becomes
\begin{equation}
\label{eq:FKPP:w(r)}
  \frac{\mathrm{d}w}{\mathrm{d}r} = {}- c + \frac{f(-r)}{\psi(w)} \,,
    \quad r\in (-1,1) \,,
\end{equation}
with the boundary conditions
\begin{equation}
\label{bc:FKPP:w(r)}
  w(-1) = w(1) = 0 \,.
\end{equation}
Function $w_i = \Psi_{-1}(z_i)$ satisfies eq.~\eqref{eq:FKPP:w(r)}
with $c_i$ in place of $c$; for $i=1,2$.
We subtract \eqref{eq:FKPP:w(r)} for $i=1$ from
\eqref{eq:FKPP:w(r)} for $i=2$, thus arriving at
\begin{equation*}
\begin{aligned}
&   \frac{\mathrm{d}}{\mathrm{d}r}\, (w_2(r) - w_1(r))
  =
{}- c_2 + \frac{f(-r)}{ \psi(w_2(r)) }
  + c_1 - \frac{f(-r)}{ \psi(w_1(r)) }
\\
& = (c_1 - c_2)
  - \frac{f(-r)}{ \psi(w_1(r))\, \psi(w_2(r)) }\,
    [ \psi(w_2(r)) - \psi(w_1(r)) ]
  \geq c_1 - c_2 \,,
\end{aligned}
\end{equation*}
thanks to $\psi(w_2(r))\geq \psi(w_1(r)) > 0$ and
$f(-r) < 0$ for every $r\in (a,b)$, i.e.,
$-r\in (-b,-a)\subset (-1,\mu)$.
Equivalently,
\begin{equation*}
  r \,\longmapsto\, (w_2(r) - w_1(r)) - (c_1 - c_2) r
  \colon [a,b]\to \RR
\end{equation*}
is a monotone increasing, continuous function.
In particular, we have
\begin{equation*}
\begin{aligned}
  (w_2(r) - w_1(r)) - (c_1 - c_2) (r-a)
  \geq
  (w_2(a) - w_1(a)) \geq 0
\\
    \quad\mbox{ for all }\, r\in [a,b] \,.
\end{aligned}
\end{equation*}
Since $c_1 > c_2$, this shows that
$w_2(r) > w_1(r)$ holds for all $r\in (a,b]$.
Finally, function
$\Psi\vert_{\RR_+}\colon \RR_+\to \RR$
being continuous and strictly monotone increasing, we conclude that
$z_2(r) = \Psi(w_2(r)) > z_1(r) = \Psi(w_1(r))$
holds for all $r\in (a,b]$, as claimed.

Our corollary is proved.
\qed
\par\vskip 10pt

Now let us return to eq.~\eqref{gen:FKPP:y(r)} with
$H = \psi\circ \Psi_{-1}$ on $\RR_+$, that is,
to our original equation, eq.~\eqref{eq:FKPP:y(r)},
\begin{equation*}
    \frac{\mathrm{d}y}{\mathrm{d}r}
  - c\, \psi\left( \Psi_{-1}(y^{+}) \right) + f(r) = 0 \,,
    \quad r\in (-1,1) \,.
\end{equation*}

The {\it uniqueness\/} of the critical speed $c\in \RR$,
for which eq.~\eqref{eq:FKPP:y(r)} possesses a positive solution
$y\colon (-1,1)\to \RR$ satisfying the Dirichlet boundary conditions
$\lim_{r\searrow -1} y(r) = 0$ and
$\lim_{r\nearrow +1} y(r) = 0$, follows from the following proposition
stated and proved for $z(r) = y(-r)$, $r\in [-1,1]$.

\begin{proposition}\label{prop-unique_c}
Let\/
$H = \psi\circ \Psi_{-1}$ on $\RR_+$ and let\/
$c_1, c_2\in \RR$ satisfy\/ $c_1\geq c_2$ and\/ $c_1\geq 0$.
Assume that\/
$z_1, z_2\colon [-1,1]\to \RR$ are two absolutely continuous functions
satisfying the differential equations
\eqref{eq:FKPP:z_1} and \eqref{eq:FKPP:z_2}, respectively,
such that\/ $z_1(r) > 0$ for all\/ $r\in (-1,1)$, together with\/
$z_i(-1) = z_i(1) = 0$ for\/ $i=1,2$.
Then we must have $c_1 = c_2$ and\/
$z_1\equiv z_2$ in $[-1,1]$.
\end{proposition}

\par\vskip 10pt
\proof
On the contrary, suppose that $c_1 > c_2$ is possible.
Corollary~\ref{cor-sub<sup}, Part~{\rm (i)}, with $a = -1$ and $b = 1$,
implies that $z_1(r) < z_2(r)$ for every $r\in (-1,1)$.
Part~{\rm (iii)} forces also $z_1(1) < z_2(1)$,
a contradiction with our boundary conditions
$z_1(1) = z_2(1) = 0$.

We have proved $c_1 = c_2$.
The equality $z_1\equiv z_2$ in $[-1,1]$
follows from Corollary~\ref{cor-sub-sup} with $c_1 = c_2$.
\qed
\par\vskip 10pt

Finally, the following lemma excludes the case $c^{\ast} = +\infty$
in eq.~\eqref{def:c^*}.
Consequently, the {\it existence\/} of the critical speed $c\in \RR$
follows from eq.~\eqref{def:c^*}
and remarks thereafter, $c = c^{\ast}\in (0,\infty)$.

\begin{lemma}\label{lem-z_c(1):-infty}
We have\/ $\lim_{c\to +\infty} z_c(1) < 0$.
\end{lemma}

\par\vskip 10pt
\proof
By Corollary~\ref{cor-sub<sup}, Part~{\rm (i)}, with $a = -1$ and $b = 1$,
the monotone decreasing limit
$L = \lim_{c\to +\infty} z_c(1)$ exists, i.e.,
$z_c(1)\searrow L$ as $c\nearrow +\infty$, and satisfies
$-\infty\leq L\leq z_0(1) < \infty$.
On the contrary, suppose that $L\geq 0$.
This forces $z_c(1)\geq L\geq 0$ for every $c\geq 0$.

Given any $c>0$, we must have $z_c(r) > 0$ for every $r\in (-1,-\mu)$.
To verify this claim, we first show that, for every
$\delta\in (0, 1-\mu)$ there is some $r_{\delta}\in (-1, -1+\delta)$,
such that $z_c(r_{\delta}) > 0$.
Namely, otherwise we would have $z_c(r)\leq 0$ for every
$r\in (-1, -1+\delta)$ and, consequently,
$z_c'(r)\equiv \frac{\mathrm{d}}{\mathrm{d}r} z_c(r) = f(-r) > 0$,
by eq.~\eqref{gen:FKPP:z(r)}, which in turn yields
$z_c(r) > 0$ for every $r\in (-1, -1+\delta)$, a contradiction.
Again, notice that
$r\in (-1, -1+\delta)$ means
$-r\in (1-\delta,1)\subset (\mu,1)$ whence $f(-r) > 0$.
Next, we show that $z_c(r) > 0$ for every
$r\in [r_{\delta}, -\mu)$.
Indeed, if $z_c(r') = 0$ for some $r'\in (r_{\delta}, -\mu)$,
then there is another number $r''\in (r_{\delta}, r']$, such that
$z_c(r'') = 0$ and $z_c(r) > 0$ for every $r\in [r_{\delta}, r'')$.
But this forces $z_c'(r'')\leq 0$ which contradicts
$z_c'(r'') = f(-r'') > 0$, by eq.~\eqref{gen:FKPP:z(r)}.
Finally, letting $\delta\searrow 0$ we arrive at
$z_c(r) > 0$ for every $r\in (-1,-\mu)$.

Since $c>0$ and $f(r) < 0$ whenever $r\in (-1,\mu)$,
it follows from eq.~\eqref{gen:FKPP:z(r)} that
$z_c'(r)\leq f(-r) < 0$ for all $r\in (-\mu,1)$.
Consequently, $z_c\colon [-1,1]\to \RR$
is strictly monotone decreasing on the interval $[-\mu,1]$.
Recalling $z_c(1)\geq L\geq 0$ and
$z_c(r) > 0$ for every $r\in (-1,-\mu)$, we conclude that
$z_c(r) > 0$ for every $r\in (-1,1)$.
Now the positivity of $z_c$ on the interval $(-1,1)$ enables us to invoke
eq.~\eqref{eq:FKPP:w(r)} for the unknown function
$w\equiv w_c = \Psi_{-1}(z_c)$ on $(-1,1)$
with the initial condition $w(-1) = 0$,
\begin{equation*}
  \frac{\mathrm{d}w}{\mathrm{d}r} = {}- c + \frac{f(-r)}{\psi(w)} \,,
    \quad r\in (-1,1) \,.
\end{equation*}
Here, $w_c(r) = \Psi_{-1}(z_c(r)) > 0$ for all $r\in (-1,1)$.
In particular, we have
$w_c'(r) < - c$ for every $r\in (-\mu,1)$.
Taking $c > 0$ sufficiently large, say, $c\geq c_0 > 0$,
we conclude that
\begin{equation*}
\begin{aligned}
    w_c(-\mu + s)
  = w_c(-\mu) + \int_{-\mu}^{-\mu + s} w_c'(r) \,\mathrm{d}r
  < w_c(-\mu) - cs\leq w_c(-\mu) - c_0 s\leq 0
\\
  \quad\mbox{ whenever }\quad
  \frac{w_c(-\mu)}{c_0}\leq s < 1 + \mu \,.
\end{aligned}
\end{equation*}
But this contradicts the fact that $w_c(r) > 0$ for all $r\in (-1,1)$.

We have proved $L<0$ as desired.
\qed
\par\vskip 10pt

Proposition~\ref{prop-TW} follows directly from a combination of
Lemma~\ref{lem-z_c(1):-infty} (existence) and
Proposition~\ref{prop-unique_c} (uniqueness).
In addition, Lemma~\ref{lem-TW} is a consequence of
the continuous dependence of the travelling wave $q$ and the speed $c$
upon the given data,
combined with a standard compactness argument
(Arzel\`a\--Ascoli's theorem),
by a continuity (convergence) result from
{\sc Ph.~Hartman} \cite[Theorem 2.1, p.~94]{Hartman}
for problem \eqref{i12}, \eqref{i13}.

\section{Convergence to the limit problem}
\label{s:conv-TW}

In analogy with the approach in {\sc E.\ Feireisl} \cite{EF14},
our proof of Theorem~\ref{thm-sol} is based on a comparison principle.
To begin with, we introduce an approximation family
$\{ \Phi_{\alpha} \}_{\alpha > 0}$
satisfying \eqref{ap1} -- \eqref{ap4}.
We start with a simple result for $N=1$.

\begin{lemma}\label{lem-approx-TW}
Suppose that\/ $v_{\alpha}$ is a weak solution of the Cauchy problem
\eqref{i1}, \eqref{i2} with $N=1$ and $\Phi = \Phi_{\alpha}$
starting from (smooth) initial data $v_0$,
\begin{align*}
  \partial_x v_0 \leq 0,\ v_0(x) = \lambda_1\
    \mbox{ for all }\ x < a,\ v_0(x) = \lambda_2
\\
    \quad\mbox{ for all }\, x > b,\ a < b,\
    \lambda_1\in (\mu,1],\ \lambda_2\in [-1,\mu) \,.
\end{align*}
Then
\begin{equation}
\label{aaa}
  \lim_{t\to \infty} v_{\alpha}(t, \underline{c} t) = 1\
    \mbox{ for any }\ \underline{c} < c ,\quad
  \lim_{t\to \infty} v_{\alpha}(t, \overline{c} t) = - 1\
    \mbox{ for any }\ \overline{c} > c ,
\end{equation}
uniformly for $\alpha\searrow 0$.
\end{lemma}
\par\vskip 10pt

\begin{remark}\label{rem-approx-TW}\nopagebreak
\begingroup\rm
Uniformly in \eqref{aaa} means that there exists
$\alpha_0 = \alpha_0 (\underline{c}, \overline{c}) > 0$
such that for any $\eps > 0$ there exist $T(\eps) > 0$ such that
\[
  v_{\alpha} (t, \underline{c} t) > 1 - \eps ,\
  v_{\alpha} (t, \overline{c} t) < -1 + \eps \
    \mbox{for all}\ \alpha < \alpha_0(\underline{c}, \overline{c}),\
    t > T(\eps).
\]
\endgroup
\end{remark}
\par\vskip 10pt

\bProof
In view of the symmetry of the problem with respect to the change
$v\approx -v$, it is enough to show
\bFormula{c1-}
\lim_{t \to \infty} v_\alpha ( t, \underline{c}t) = 1 \ \mbox{for any}\ \underline{c} < c, \ \mbox{uniformly for} \ \alpha \searrow 0,
\eF
for
\[
\partial_x v_0 \leq 0,\
v_0 = v_0 (x) = \lambda_1 \ \mbox{for all}\ x < a,\ v_0(x) = -1  \ \mbox{for all}\ x > b,\
a < b, \ \lambda_1 \in (\mu, 1).
\]

\medskip

{\bf Step~1.}$\;$
We show that for any $\eps > 0$ there exist a time $T(\eps)$, $a(\eps)$, and $\alpha_0 > 0$ such that
\bFormula{c1}
v_\alpha (T(\eps), x) > 1 - \eps \ \mbox{for all}\ x < a(\eps), \ \alpha < \alpha_0.
\eF

By the comparison principle, the spatially homogeneous solution $\overline{v} = \overline{v}(t)$,
\[
\partial_t \overline{v}(t) = f(\overline{v}(t)), \ \overline{v}(0) = \lambda_1
\]
dominates $v_\alpha $,
\[
\overline{v}(t) \geq v_\alpha (x,t) \ \mbox{for all}\ t \geq 0,\ x \in R, \ \alpha > 0.
\]
Since $f$ satisfies (\ref{i3}), we have
\begin{equation*}
\label{c2}
  \overline{v}(t)\to 1 \quad\mbox{ as }\ t\to \infty \,.
\end{equation*}

Given parameters $\delta\in (0,\infty)$ and $Y\in \mathbb{R}$,
let us consider an auxiliary function
$\omega_{\delta,Y}\colon \mathbb{R}\to [0,1]$ defined for every
$x\in \mathbb{R}$ by
$\omega_{\delta,Y}(x)\eqdef \omega_{\delta,0}(x-Y)$
where
\begin{equation*}
  \omega_{\delta,0}(x)\eqdef
  \left\{
\begin{alignedat}{2}
  & \qquad 1
  & \quad\hbox{ if }\;
  & |x|\leq \textstyle{ \frac{1}{2} } \,;
\\
  & \delta\left( \textstyle{ \frac{1}{2} + \frac{1}{\delta} } - |x|
          \right)
  & \quad\hbox{ if }\;
  & \textstyle{ \frac{1}{2} } < |x|\leq
    \textstyle{ \frac{1}{2} + \frac{1}{\delta} } \,;
\\
  & \qquad 0
  & \quad\hbox{ if }\;
  & \textstyle{ \frac{1}{2} + \frac{1}{\delta} }
    < |x| < \infty \,.
\end{alignedat}
\right.
\end{equation*}
Clearly, $\omega_{\delta,0}$ is an even function (about $0$), i.e.,
$\omega_{\delta,0}(x) = \omega_{\delta,0}(|x|)$
for all $x\in \mathbb{R}$.

Using equation \eqref{i1} with $\Phi$ replaced by $\Phi_\alpha$, we get
\bFormula{c2a}
\int_R \omega_{\delta,Y} (\overline{v} - v_\alpha) (\tau) \ \dx
\eF
\[
\leq  \int_R \omega_{\delta,Y} (\overline{v}(0) - v_0 ) \ \dx +  \int_0^\tau \int_R \omega_{\delta,Y} \left| f(\overline{v}) -
f(v_\alpha) \right| \,\mathrm{d}x \,\mathrm{d}t + \int_0^\tau \int_R \partial_z \varphi_\alpha (|\partial_x v_\alpha |)|\partial_x \omega_{\delta,Y} |
  \,\mathrm{d}x \,\mathrm{d}t
\]
for any $\tau \geq 0$.

In accordance with Corollary~\ref{cor-int_reg},
there is $M$, $\alpha_0 > 0$ such that
\begin{equation*}
  |\partial_x v_\alpha(t,x) | \leq M
    \quad\mbox{ for all }\
    t\geq 0,\ x\in \mathbb{R}^N ,\ \alpha < \alpha_0 \,.
\end{equation*}
Moreover,
by virtue of the hypothesis \eqref{p3}, we get
\begin{equation*}
\begin{aligned}
       \mathrm{d}_z \log (z^{\Lambda_1})
     = \frac{\Lambda_1}{z}
  \leq \mathrm{d}_z \left( \log \mathrm{d}_z \varphi(z) \right)
  \leq \frac{\Lambda_2}{z}
     = \mathrm{d}_z \log (z^{\Lambda_2})
\\
  \quad\mbox{ for all }\, z\in (0,\infty)\,,
\end{aligned}
\end{equation*}
which upon integration shows that, with some constant $m\in (0,1)$,
on every compact interval $[0,M]\subset \RR_+$ we have
\begin{equation*}
  0\leq \partial_z\varphi_\alpha(z)
  \leq \mathrm{const}(M) \left( z + z^m \right)
    \quad\mbox{ for all }\
    0\leq z\leq M,\ \alpha < \alpha_0 \,;
\end{equation*}
hence
\begin{equation*}
  \left\| \partial_z \varphi_\alpha(|\partial_x v_\alpha|)
  \right\|_{L^1 + L^q(\RR)}
  \leq \mathrm{const}(M)
    \left( \| \partial_x v_\alpha\|_{L^1(\RR)}
  + \| \partial_x v_\alpha \|_{L^1(\RR)}^{m} \right) \,,
    \quad q = \frac{1}{m} \,.
\end{equation*}
Since $\partial_x v_\alpha \leq 0$, we get
\begin{equation*}
  \| \partial_x v_\alpha(t,\,\cdot\,) \|_{L^1(\RR)} \leq 2
  \quad\mbox{ for all }\ t\geq 0,\ \alpha < \alpha_0 \,,
\end{equation*}
thus yielding
\begin{equation*}
  \left\| \partial_z \varphi_\alpha(|\partial_x v_\alpha|)
  \right\|_{L^1 + L^q(\RR)}\leq \mathrm{const}(M) ,\ q = \frac{1}{m} \,.
\end{equation*}
We remark that the sum
$L^1 + L^q(\RR)\equiv L^1(\RR) + L^q(\RR)$
and its norm are defined in a standard way used
in interpolation theory, cf.\
{\sc H.\ Triebel} \cite[{\S}1.2.1]{Triebel}.

On the other hand,
\[
  \| \partial_x \omega_{\delta, Y} \|_{L^\infty(R)} = \delta \,,\quad
  \| \partial_x \omega_{\delta, Y} \|_{L^1(R)} = 2 \,;
\]
therefore, going back to (\ref{c2a}), we may infer that
\bFormula{c4}
\int_R \omega_{\delta, Y} (\overline{v} - v_\alpha ) (\tau) \ \dx \leq \exp( L_f \tau)
\left( \int_R \omega_{\delta, Y} (\overline{v}(0) - v_0 ) \ \dx + \tau \chi(\delta) \right),
\eF
where $L_f$ is the Lipschitz constant of $f$ and $\chi(\delta) \to 0$ for $\delta \searrow 0$.

In accordance with \eqref{c2},
we may fix $\tau = T(\eps)$ in (\ref{c4}) so that
\[
\overline{v}(t) > \left( 1 - \frac{\eps}{2} \right) \ \mbox{for all}\ t > \tau,
\]
and take $\delta > 0$ so small and $Y < 0$ so that (\ref{c4}) yields the existence of a point
$a(\eps)$ such that
\[
v_\alpha (T(\eps), a(\eps)) > 1 - \eps, \ \alpha < \alpha_0.
\]
As $\partial_x v_\alpha \leq 0$, the desired conclusion (\ref{c1}) follows.

\medskip

{\bf Step~2.}$\;$
In agreement with the previous discussion, it is enough to examine the initial datum
\[
\partial_x v_0 \leq 0,\
v_0 (x) = 1 - \eps \ \mbox{for all}\ x < a,\ v_0(x) = -1  \ \mbox{for all}\ x > b,\
a < b,
\]
where $\eps \searrow 0$.

Keeping the conclusion of Lemma \ref{lem-TW} in mind, we take a family $f_\eps \approx f$ satisfying the hypotheses of Lemma \ref{lem-TW} with
\[
b_\eps = 1 - \eps, \ a_\eps = -1 - \eps, \ f_\eps \leq f.
\]
Consequently, taking $\alpha_0 (\underline{c}) > 0$ small enough, we can find a traveling wave $q_\alpha$ with the propagation speed $c_\alpha$,
\begin{equation*}
  c > c_\alpha > \underline{c}
  \quad\mbox{ for all }\ \alpha < \alpha_0(\underline{c}) \,,
\end{equation*}
such that
\[
q_\alpha ( x + D) \leq v_0(x) \ \mbox{for all}\ x \in R
\]
for a suitable constant $D \in R$, where, by comparison,
\[
v_\alpha (t,x) \geq q_{\alpha }( x + D - c_\alpha t ), \ \mbox{in particular,} \ v_\alpha
(t, \underline{c}t) \geq q_{\alpha}(  D + (\underline{c} - c_\alpha) t ) \to 1 - \eps
\ \mbox{as} \ t \to \infty \
\]
uniformly for $\alpha < \alpha_0$.

Thus we have shown
\[
\liminf_{t \to \infty} v_\alpha (t, \underline{c}t) \geq 1 - \eps \ \mbox{uniformly for}\ \alpha < \alpha_0.
\]
Since $\eps > 0$ can be taken arbitrarily small, the desired conclusion follows.

\qed

The next step is to extend the previous result to the case of
radially symmetric data in $\RR^N$.
To this end, we introduce a new variable $r = |x|$ and rewrite (formally)
equation \eqref{i1} for the radially symmetric solutions:
\begin{align}
\label{c5}
\partial_t u = \partial_r \left( \partial \varphi (|\partial_r u|) \frac{ \partial_r u }{|\partial_r u|} \right) +
\frac{N-1}{r} \partial \varphi (|\partial_r u|) \frac{ \partial_r u }{|\partial_r u|}   + f(u),
\\
\nonumber
  u = u(t,r), \ r > 0, \ \partial_r u(t,0) = 0.
\end{align}
Here, we have used the following simple relations for
the radially symmetric functions
$\Phi(Z)\equiv \varphi(|Z|)$, $Z\in \RR^N$, and
$u(x,t)\equiv u(r,t)$, $x\in \RR^N$, $r = |x|$:
\begin{align*}
& \partial\Phi(Z) = \partial\varphi(|Z|) \,\frac{Z}{|Z|}
    \quad\mbox{ and }\quad
  \nabla_x u(x,t) = \partial_r u(r,t) \,\frac{x}{r}
    \quad\mbox{ for }\, Z,x\in \RR^N\setminus \{\mathbf{0}\} \,,
\\
& \nabla_x (|x|) = \frac{x}{|x|} \quad\mbox{ and }\quad
  \Div \genfrac{(}{)}{}0{x}{|x|} = \frac{N-1}{|x|}
    \quad\mbox{ for }\, x\in \RR^N\setminus \{\mathbf{0}\} \,.
\end{align*}

\begin{lemma}\label{lem-c2}
Suppose that $v_\alpha = v_\alpha (t,r)$ is the radially symmetric weak solution of the Cauchy problem \eqref{i1}, \eqref{i2}, with
$\Phi_\alpha$ emanating from
(smooth) initial datum $v_0 = v_0(r)$,
\begin{align*}
& \partial_r v_0\leq 0,\
  v_0 (r) = \lambda \quad\mbox{ for all }\ r\in (0,R) \,,\quad
\\
& v_0(r)\geq -1 \quad\mbox{ for all }\ r\in (\tilde{R},\infty) ,\
  0 < R < \tilde{R} \,,
    \quad\mbox{  with }\ \lambda\in (\lambda_0,1),\
    \lambda_0\in (\mu,1) \,.
\end{align*}

Then for any $\eps > 0$ and $\underline{c} < c$, there exists $\alpha = \alpha_0 (\underline{c})$, a time $T = T(\lambda_0, \eps, \underline{c})$ and $R_0 = R_0(\lambda_0, \eps, \underline{c})$ such that
\[
v_\alpha (t,r) \geq 1 - \eps \ \mbox{for all} \ t \in [T, 2T],\ |x| < R + \underline{c} t ,\ \alpha < \alpha_0,
\]
whenever $R > R_0$.
\end{lemma}

\bProof

The proof is along the same lines as that of Lemma \ref{lem-approx-TW}.
Assuming $\delta > 0$ is chosen small enough, we fix a (smooth) profile
\[
w_0 = w_0(x), \ \partial_x w_0 \leq 0, \ w_0(x) = \lambda \ \mbox{for}\ x < 0,\ w_0(x) = -1  \ \mbox{for} \ x > 1,
\]
and such that
\[
w_0(r - R) \leq v_0(r) \ \mbox{for all}\ r > 0.
\]
By virtue of the comparison principle, it is enough to show the conclusion of the lemma for $v_0(r) = w_0(r-R)$.

\medskip

{\bf Step~1.}$\;$
Consider the unique solution $w_\alpha$ of the Cauchy problem
\[
\partial_t w_\alpha = \partial_x \left( \partial_z \varphi_\alpha (|\partial_x w_\alpha |) \frac{\partial_x w_\alpha }{|\partial_x w_\alpha |} \right) + f(w_\alpha),\ t > 0,\ x \in R^1, w_\alpha(0,x) = w_0(x).
\]
Making use once more of the comparison principle we deduce that
\[
w_\alpha (t,r - R) \geq v_\alpha (t,r) \ \mbox{for any}\ r > 0 \ \mbox{and}\ t \geq 0,
\]
where, by virtue of Lemma \ref{lem-approx-TW},
\bFormula{gr1}
\lim_{t \to \infty} \inf_{r < \underline{c} t} w_\alpha (t,r) = 1 \ \mbox{for any}\ \underline{c} < c
\ \mbox{uniformly in}\ \alpha < \alpha_0 (\underline{c}).
\eF

\medskip

{\bf Step~2.}$\;$
Similarly to the proof of Lemma \ref{lem-approx-TW},
we take the function $\omega_{\delta,Y}$ and compute the ``distance''
\[
\int_0^\infty \omega_{\delta, Y} \Big[ w_\alpha (\tau, r- R) - v_\alpha (\tau,r) \Big] \ {\rm d}r \leq \int_0^\tau \int_0^\infty \omega_{\delta, Y} \Big| f(w_\alpha (t, r-R)) - f(v_\alpha (t,r)) \Big| \ {\rm d}r \ {\rm d}t
\]
\[
+ \int_0^\tau \int_0^\infty \Big( \partial_z \varphi_\alpha (|\partial_r v_\alpha (t,r)|) + \partial_z \varphi_\alpha (|\partial_r
w_\alpha (t,r - R)|) \Big) |\partial_r \omega_{\delta, Y}|
\ {\rm d}r \ {\rm d}t
\]
\[
 + \int_0^\tau \int_0^\infty \frac{N-1}{r} \omega_{\delta,Y}  \partial_z \varphi_\alpha (|\partial_r v_\alpha (t,r)|) \ {\rm d}r \ {\rm d}t
\]
for $Y > 1 + \frac{1}{\delta}$.

Now, exactly as in the proof of Lemma \ref{lem-approx-TW},
we may show that
\[
\left| \int_0^\tau \int_0^\infty \Big( \partial_z \varphi_\alpha (|\partial_r v_\eps (t,r)|) + \partial_z \varphi_\alpha (|\partial_r w_\eps (t,r - R)|) \Big) |\partial_r \omega_{\delta, Y}|
\ {\rm d}r \ {\rm d}t \right| \leq \tau \chi(\delta), \ \mbox{with}\ \chi(\delta) \to 0\ \mbox{as}\ \delta \searrow 0,
\]
uniformly for $\alpha < \alpha_0$.

Consequently, by Gronwall's lemma,
\bFormula{gr2}
\int_0^\infty \omega_{\delta, Y} \Big[ w_\alpha (\tau, r- R) - v_\alpha (\tau,r) \Big] \ {\rm d}r \leq \exp \left( L_f \tau \right) \left(
\tau \chi(\delta) + \int_0^\tau \int_0^\infty \frac{N-1}{r} \omega_{\delta,Y}  \partial_z \varphi_\alpha (|\partial_r v_\alpha (t,r)|) \ {\rm d}r \ {\rm d}t \right)
\eF
for $\tau > 0$, $Y > 1 + \frac{1}{\delta}$.

\medskip

{\bf Step~3.}$\;$
In accordance with (\ref{gr1}), there exists $T = T(\lambda_0, \eps, \underline{c})$ such that
\bFormula{gr3}
w_\alpha (t, r - R) \geq 1 - \frac{\eps}{4} \ \mbox{for all}\ t > T \ \mbox{and} \ r < R + \underline{c} t + 2,\ \alpha < \alpha(\underline{c}).
\eF

Next, fix $\delta > 0$ so that
\bFormula{gr4}
2 \chi(\delta) T \exp (2 L_f T) \leq \frac{\eps}{4}.
\eF

Finally, take $Y \geq Y_0(\delta)$ so large that
\bFormula{gr5}
\int_0^{2T} \int_0^\infty \frac{N-1}{r} \omega_{\delta,Y}  \partial_z \varphi_\alpha (|\partial_r v_\alpha(t,r)|) \ {\rm d}r \ {\rm d}t \leq \frac{\eps}{4},\ Y \geq Y_0.
\eF

Thus, for $R \geq R_0$ large enough so that
\[
Y = R + \underline{c} \tau + 1 > Y_0 \ \mbox{for all}\ \tau \in [T,2T]
\]
in the inequality (\ref{gr2}), we may use (\ref{gr3} - \ref{gr5}), together with the monotonicity of $v_\alpha $ in $r$, to obtain the desired conclusion.

\qed

Applying Lemma~\ref{lem-c2} recursively on the sequence of time intervals $[nT,(n+2)T]$ we obtain:

\bCorollary{c1}
Suppose that $v_\alpha = v_\alpha(t,r)$ is the radially symmetric weak solution of the Cauchy problem (\ref{i1}), (\ref{i2}),
with $\Phi_\alpha$, emanating from
(smooth) initial datum $v_0 = v_0(r)$,
\begin{align*}
& \partial_r v_0\leq 0,\
  v_0 (r) = \lambda \quad\mbox{ for all }\ r\in (0,R) \,,\quad
\\
& v_0(r)\geq -1 \quad\mbox{ for all }\ r\in (\tilde{R},\infty) ,\
  0 < R < \tilde{R} \,,
    \quad\mbox{  with }\ \lambda\in (\lambda_0,1),\
    \lambda_0\in (\mu,1) \,.
\end{align*}

Then for any $\eps > 0$ and $\underline{c} < c$, there exist $\alpha_0 = \alpha_0 (\underline{c})$, a time $T = T(\lambda_0, \eps, \underline{c})$ and $R_0 = R_0(\lambda_0, \eps, \underline{c})$ such that
\[
v_\alpha (\tau,r) \geq 1 - \eps \ \mbox{for all} \ \tau > T ,\ 0< r < R + \underline{c} \tau,\ \alpha < \alpha_0,
\]
as long as
\[
R + \underline{c} t > R_0 \ \mbox{for all}\ t \in [0, \tau].
\]
\eC

Introducing $u_\eps = u \left( \frac{t}{\eps}, \frac{x}{\eps} \right)$, the function $u_\eps$ solves the scaled equation \eqref{i4}.
Adapting Corollary~\ref{Cc1} we get:

\bCorollary{c2}
Suppose that $\{ u_{\alpha,\eps } \}_{\eps > 0}$ is a family of solution of the scaled equation \eqref{i4}, with $\Phi = \Phi_\eps$,
emanating from the initial data
\[
-1 \leq u_{\eps, \alpha}(0, \cdot) \leq 1, \ u_{\eps, \alpha}(0, x) \geq \lambda \in (\mu, 1] \ \mbox{for all}\ |x| < R.
\]

Then, given a compact set
\[
\mathcal{K} \subset \Big\{ t > 0, \ |x| < R + ct \Big\}
\]
there is $\alpha_0 (\mathcal{K}) > 0$ such that
\[
\lim_{\eps \to 0} u_{\eps, \alpha} (t,x) = 1 \ \mbox{uniformly in} \ \mathcal{K} \ \mbox{and uniformly for}\ \alpha <
\alpha_0 (\mathcal{K}).
\]
\eC

Finally, using ``symmetric'' arguments we get

\bCorollary{c3}
Suppose that $\{ u_{\alpha,\eps } \}_{\eps > 0}$ is a family of solution of the scaled equation (\ref{i4}), with $\Phi = \Phi_\eps$, emanating from the initial data
\[
-1 \leq u_{\eps, \alpha}(0, \cdot) \leq 1, \ u_{\eps, \alpha}(0, x) \leq \lambda \in [-1, \mu) \ \mbox{for all}\ |x| < R.
\]

Then, given a compact set
\[
\mathcal{K} \subset \Big\{ t > 0, \ |x| < R - ct \Big\}
\]
there is $\alpha_0 (\mathcal{K}) > 0$ such that
\[
\lim_{\eps \to 0} u_{\eps, \alpha} (t,x) = -1 \ \mbox{uniformly in} \ \mathcal{K} \ \mbox{and uniformly for}\ \alpha <
\alpha_0 (\mathcal{K}).
\]

\eC

The final observation is that
Corollaries \ref{Cc2}, \ref{Cc3}
imply the conclusion of Theorem~\ref{thm-sol}.
The proof is exactly the same as in \cite[Section 5]{EF14}
provided we consider only those weak solutions
that have been suitably introduced in Definition~\ref{def-sol}.

\bigskip
\bigskip
\bigskip



%
%
\makeatletter \renewcommand{\@biblabel}[1]{\hfill#1.} \makeatother
%
%

\end{document}